\documentclass[12pt]{amsart}
\usepackage{amsfonts}

\usepackage{amscd}
\newtheorem{theorem}{Theorem}[section]

\newtheorem{lemma}[theorem]{Lemma}

\newtheorem{remark}[theorem]{Remark}

\newtheorem{definition}[theorem]{Definition}

\def\limfunc#1{\mathop{\rm #1}}%
\def\func#1{\mathop{\rm #1}\nolimits}%
\long\def\TeXButton#1#2{#2}%

\calclayout
\numberwithin{equation}{section}
\allowdisplaybreaks
\theoremstyle{definition}
\theoremstyle{remark}
\numberwithin{equation}{section}

\begin{document}
\title[Dirichlet Problem]{The Dirichlet problem for elliptic equations in divergence and nondivergence
form with singular drift term}
\author{Cristian Rios}
\address{Trinity College}
\email{Cristian.Rios@trincoll.edu}
\date{
\today%
}
\maketitle

\begin{abstract}
Given two elliptic operators $\mathcal{L}_{0}$ and $\mathcal{L}_{1}$ in
nondivergence form, with coefficients $\mathbf{A}_{\ell }$ and drift terms $%
\mathbf{b}_{\ell }$, $\ell =0,1$ satisfying 
\begin{equation*}
\sup_{\left| Y-X\right| \le \frac{\delta \left( X\right) }{2}}\frac{\left| 
\mathbf{A}_{0}\left( Y\right) -\mathbf{A}_{1}\left( Y\right) \right|
^{2}+\delta \left( X\right) ^{2}\left| \mathbf{b}_{0}\left( Y\right) -%
\mathbf{b}_{1}\left( Y\right) \right| ^{2}}{\delta \left( X\right) }dX
\end{equation*}
is a Carleson measure in a Lipschitz domain $\Omega \subset \Bbb{R}^{n+1}$, $%
n\ge 1$, (here $\delta \left( X\right) =\limfunc{dist}\left( X,\partial
\Omega \right) $). If the harmonic measure $d\omega _{\mathcal{L}_{0}}\in
A_{\infty }$, then $d\omega _{\mathcal{L}_{1}}\in A_{\infty }$. This is an
analog to Theorem 2.17 in \cite{HL} for divergence form operators. \newline
As an application of this, a new approximation argument and known results we
obtain: Let $\mathcal{L}$ be an elliptic operator with coefficients $\mathbf{%
A}$ and drift term $\mathbf{b}$; $\mathcal{L}$ can be in divergence or
nondivergence form. If 
\begin{equation*}
\sup_{\left| Y-X\right| \le \frac{\delta \left( X\right) }{2},\left|
Z-X\right| \le \frac{\delta \left( X\right) }{2}}\frac{\left| \mathbf{A}%
\left( Y\right) -\mathbf{A}\left( Z\right) \right| ^{2}+\delta \left(
X\right) ^{2}\left| \mathbf{b}\left( Y\right) -\mathbf{b}\left( Z\right)
\right| ^{2}}{\delta \left( X\right) }dX
\end{equation*}
is a Carleson measure in $\Omega $, then $d\omega _{\mathcal{L}}\in
A_{\infty }$.This extends the results in \cite{KP} for divergence form
operators while provides totally new results for nondivergence form
operators. The results are sharp in all cases.
\end{abstract}

\section{Introduction and Background}

Given a bounded Lipschitz domain $\Omega \subset \Bbb{R}^{n+1}$, $n\ge 1$,
and an operator $\mathcal{L}$ given by 
\begin{eqnarray*}
\mathcal{L} &=&\mathcal{L}_{d}=\limfunc{div}\mathbf{A}\nabla ,\qquad \text{%
(divergence fom)\qquad or} \\
\mathcal{L} &=&\mathcal{L}_{n}=\mathbf{A}\cdot \nabla ^{2},\qquad \text{%
(nondivergence form),}
\end{eqnarray*}
the harmonic measure at $X\in \Omega $, $d\omega _{\mathcal{L}}^{X}$, is the
unique Borel measure on $\partial \Omega $ such that for all continuous
functions $g\in \mathcal{C}\left( \partial \Omega \right) $, 
\begin{equation*}
u\left( X\right) =\int_{\partial \Omega }g\left( Q\right) d\omega _{\mathcal{%
L}}^{X}\left( Q\right)
\end{equation*}
is continuous in $\overline{\Omega }$ and it is the unique solution to the
Dirichlet problem 
\begin{equation}
\left\{ 
\begin{array}{ll}
\mathcal{L}u=0 & \text{in }\Omega \\ 
u=g & \text{on }\partial \Omega .
\end{array}
\right.  \label{DP}
\end{equation}
Where we assume that the equality $\mathcal{L}u=0$ holds in the weak sense
for divergence form operators and in the strong $a.e.$ sense for
nondivergence form operators. Here $\mathbf{A}=\mathbf{A}\left( X\right) $
is a symmetric $\left( n+1\right) \times \left( n+1\right) $ matrix with
bounded measurable entries, satisfying a uniform ellipticity condition 
\begin{equation}
\lambda \left| \xi \right| ^{2}\le \xi \cdot \mathbf{A}\left( X\right) \xi
\le \Lambda \left| \xi \right| ^{2},\qquad X,\xi \in \Bbb{R}^{n+1},
\label{ellipticity}
\end{equation}
for some positive constants $\lambda $, $\Lambda $. In the nondivergence
case the entries of the matrix $\mathbf{A}$ are assumed to belong to $%
\limfunc{BMO}\nolimits\left( {\Omega }\right) $ with small enough norm. For
a given operator $\mathcal{L}$, the harmonic measures $d\omega _{\mathcal{L}%
}^{X}$, $X\in \Omega $, are regular probability measures which are mutually
absolutely continuous with respect to each other. That is, 
\begin{equation*}
k\left( X,Y,Q\right) =\frac{d\omega _{\mathcal{L}}^{Y}}{d\omega _{\mathcal{L}%
}^{X}}\left( Q\right) \in L^{1}\left( d\omega _{\mathcal{L}}^{X},\partial
\Omega \right) ,\qquad X,Y\in \Omega ,\quad Q\in \partial \Omega .
\end{equation*}
By the Harnack's principle the kernel function $k\left( X,Y,Q\right) $ is
positive and uniformly bounded in compact subsets of $\Omega \times \Omega
\times \partial \Omega $. As a consequence, to study differentiability
properties of the family $\left\{ d\omega _{\mathcal{L}}^{X}\right\} _{X\in
\Omega }$ with respect to any other Borel measure $d\nu $ on $\partial
\Omega $, it is enough to fix a point $X_{0}\in \Omega $ and study $\frac{%
d\omega }{d\nu }$, where $d\omega =d\omega _{\mathcal{L}}^{X_{0}}$ is
referred as \emph{the} harmonic measure of $\mathcal{L}$ on $\partial \Omega 
$. The well definition of the harmonic measure follows from Riesz
representation theorem if there is unique solvability of the continuous
Dirichlet problem and a boundary Maximum Principle is available.

\begin{definition}
\label{CDP}Given an elliptic operator $\mathcal{L}$, we say that the
continuous Dirichlet problem is uniquely solvable in $\Omega $, and we say
that $\mathcal{CD}$ holds for $\mathcal{L}$, if for every continuous
function $g$ on $\partial \Omega $, there exists a unique solution $u$ of (%
\ref{DP}), such that $u\in \mathcal{C}^{0}\left( \overline{\Omega }\right)
\bigcap W^{2,p}\left( \Omega \right) $ for some $1\le p\le \infty $.
\end{definition}

Given two regular Borel measures $\mu $ and $\omega $ in $\partial \Omega $, 
$d\omega \in A_{\infty }\left( d\sigma \right) $ if there exist constants $%
0<\varepsilon ,\,\delta <1$ such that for any boundary ball $\Delta =\Delta
_{r}\left( Q\right) $ and any Borel set $E\subset \Delta $, 
\begin{equation*}
\frac{\mu \left( E\right) }{\mu \left( \Delta \right) }<\delta
\Longrightarrow \frac{\omega \left( E\right) }{\omega \left( \Delta \right) }%
<\varepsilon .
\end{equation*}
The relation $d\omega \in A_{\infty }\left( d\mu \right) $ is an equivalence
relation \cite{Mu}, and any two measures related by the $A_{\infty }$
property are mutually absolutely continuous with respect to each other. From
classic theory of weights, if $d\omega \in A_{\infty }\left( d\mu \right) $
then there exists $1<q<\infty $ such that the density $h=\frac{d\omega }{%
d\mu }$ satisfies a reverse H\"{o}lder inequality with exponent $q$: 
\begin{equation*}
\left\{ \frac{1}{\mu \left( \Delta \right) }\int_{\Delta }h^{q}\,d\mu
\right\} ^{\frac{1}{q}}\le C\frac{1}{\mu \left( \Delta \right) }\int_{\Delta
}h\,d\mu .
\end{equation*}
This property is denoted $d\omega \in B_{q}\left( d\mu \right) $, and $%
d\omega \in B_{q}\left( d\mu \right) $ is equivalent to the fact that the
Dirichlet problem (\ref{DP}) for the operator $\mathcal{L}$ is solvable in $%
L^{p}\left( d\mu ,\partial \Omega \right) $, $\frac{1}{p}+\frac{1}{q}=1$
(see \cite{FKP} for details). When $\mu =\sigma $, the Euclidean measure, we
write $A_{\infty }\,$for $A_{\infty }\left( d\sigma \right) $.

\begin{definition}
\label{DPp}Let $\mathcal{L}$ be an elliptic operator that satisfies $%
\mathcal{CD}$ and let $\mu $ be a doubling measure in $\partial \Omega $. We
say that the $L^{p}\left( d\mu \right) $-Dirichlet problem is uniquely
solvable in $\Omega $, and we write that $\mathcal{D}_{p}\left( d\mu \right) 
$ holds for $\mathcal{L}$, if for every continuous function $g$ on $\partial
\Omega $, the unique solution $u$ of (\ref{DP}) satisfies 
\begin{equation*}
\left\| Nu\right\| _{L^{p}\left( d\mu \right) }\le C\left\| g\right\|
_{L^{p}\left( d\mu \right) },
\end{equation*}
for some constant independent of $g$. Here $Nu$ denotes the nontangential
maximal function of $u$ on $\partial \Omega $.
\end{definition}

In the remarkable work \cite{FKP}, the authors established a perturbation
result relating the harmonic measures of two operators in divergence form.
The analogue result was later obtained by the author in \cite{CR} for
nondivergence form operators. Before stating the result we need a few more
definitions.

\begin{definition}
\label{Carleson}Let $\Omega $ be an open set in $\Bbb{R}^{n+1}$ and let $\mu 
$ be a nonnegative Borel measure on $\partial \Omega $. For $X\in \partial
\Omega $ and $r>0$ denote by $\triangle _{r}\left( X\right) =\left\{ Z\in
\partial \Omega :\left| Z-X\right| <r\right\} $ and $T_{r}\left( X\right)
=\left\{ Z\in \Omega :\left| Z-X\right| <r\right\} $. Given a nonnegative
Borel measure $\nu $ in $\Omega $, we say that $\nu $ is a \emph{Carleson
measure} \emph{in} $\Omega $ \emph{with respect to }$\mu $, if there exist a
constant $C_{0}$ such that for all $X\in \partial \Omega $ and $r>0$, 
\begin{equation*}
\nu \left( T_{r}\left( X\right) \right) \le C_{0}\mu \left( \triangle
_{r}\left( X\right) \right) ,
\end{equation*}
The infimum of all the constants $C_{0}$ such that the above inequality
holds for all $X\in \partial \Omega $ and $r>0$ is called \emph{the Carleson
norm of }$\nu $\emph{\ with respect to }$\mu $\emph{\ in }$\Omega $. For
conciseness, we will write $\nu \in \frak{C}\left( d\mu ,\Omega \right) $
when $v$ is a Carleson measure in $\Omega $, and we denote by $\left\| \nu
\right\| _{\frak{C}\left( d\mu ,\Omega \right) }$ its Carleson norm.When $%
\mu =\sigma $ is the Lebesgue measure on $\partial \Omega $ we just say that 
$\nu $\thinspace is a Carleson norm in $\Omega $.
\end{definition}

Throughout this work, $Q_{\gamma }\left( X\right) $ denotes a cube centered
at $X$ with faces parallel to the coordinate axes and sidelength $\gamma t$;
i.e. 
\begin{equation*}
Q_{\gamma }\left( X\right) =\left\{ Y=\left( y_{1},\dots ,y_{n+1}\right) \in 
\Bbb{R}^{n+1}:\left| y_{i}-x_{i}\right| <\frac{\gamma t}{2},i=1,\dots
,n+1\right\} .
\end{equation*}

\begin{theorem}
\label{FKP0}\cite{FKP}-\cite{CR}Let $\mathcal{L}_{d,0}=\limfunc{div}\mathbf{A%
}_{0}\nabla $ and $\mathcal{L}_{d,1}=\limfunc{div}\mathbf{A}_{1}\nabla $ be
two elliptic operators with bounded measurable coefficients in $\Omega $,
and let $\omega _{d,0}$ and $\omega _{d,1}$ denote their respective harmonic
measures. Let $\sigma $ be a doubling measure on $\partial \Omega $ and
suppose that 
\begin{equation}
\delta \left( X\right) ^{-1}\sup_{Y\in Q_{\frac{\delta \left( X\right) }{2%
\sqrt{n}}}}\left| \mathbf{A}_{0}\left( Y\right) -\mathbf{A}_{1}\left(
Y\right) \right| ^{2}dX\in \frak{C}\left( d\sigma ,\Omega \right) .
\label{CarlFKP}
\end{equation}
If $d\omega _{d,0}\in A_{\infty }$ then $d\omega _{d,1}\in A_{\infty }$.
Also, if $\mathcal{CD}$ holds for the operators $\mathcal{L}_{n,0}=\mathbf{A}%
_{0}\cdot \nabla ^{2}$ and $\mathcal{L}_{n,1}=\mathbf{A}_{1}\cdot \nabla ^{2}
$ (see Definition \ref{CDP}), then their respective harmonic measures $%
\omega _{n,0},$ $\omega _{n,1}$, satisfy $d\omega _{n,0}\in A_{\infty
}\Rightarrow d\omega _{n,1}\in A_{\infty }$.
\end{theorem}

The theorems above were stated in terms of the supremum of the differences $%
\left| \mathbf{A}_{0}\left( Y\right) -\mathbf{A}_{1}\left( Y\right) \right| $
for $Y$ in Euclidean balls $\left| Y-X\right| <\frac{\delta \left( X\right) 
}{2}$. The above formulation is equivalent.

\begin{remark}
\label{exCD}By Theorem \ref{CRTth} below \cite{CRt}, a sufficient condition
for $\mathcal{CD}$ to hold for a nondivergence form operator $\mathcal{L}=%
\mathbf{A}\cdot \nabla ^{2}$ is that there exists $\rho >0$ depending on $n$
and the ellipticity constants such that 
\begin{equation}
\left\| \mathbf{A}\right\| _{\limfunc{BMO}\left( \Omega \right) }\le \rho .
\label{Arho}
\end{equation}
It is not known whether or not the continuous Dirichlet problem is uniquely
solvable in the case of elliptic nondivergence form operators with just
bounded measurable coefficients. On the other hand, even under the
restrictions (\ref{Arho}) for any $\rho >0$ it is known \cite{S} that the
continuous Dirichlet problem has non-unique ``good solutions''. That is ,
for any $\rho >0$ there exists $\mathbf{A}\left( X\right) \in $ $\func{BMO}%
\left( \Omega \right) $ with $\left\| \mathbf{A}\right\| _{\limfunc{BMO}%
\left( \Omega \right) }\le \rho $ and two sequences of $\mathcal{C}^{\infty }
$ symmetric matrices $\mathbf{A}_{0,j}$ and $\mathbf{A}_{1,j}$ with the same
ellipticity constants as $\mathbf{A}$, such that $\mathbf{A}_{\ell ,j}\left(
X\right) \rightarrow \mathbf{A}\left( X\right) $ as $j\rightarrow \infty $
for $a.e.X$, $\ell =0,1$, and such that for some continuous function $g\,$on 
$\partial \Omega $ the solutions $u_{0,j}$ and $u_{1,j}$ to the Dirichlet
problems 
\begin{equation*}
\left\{ 
\begin{array}{rrrr}
\mathcal{L}_{0,j}u_{0,j} & = & 0 & \text{in }\Omega  \\ 
u_{0,j} & = & g & \text{on }\partial \Omega .
\end{array}
\right. \quad \text{and}\quad \left\{ 
\begin{array}{rrrr}
\mathcal{L}_{1,j}u_{1,j} & = & 0 & \text{in }\Omega  \\ 
u_{1,j} & = & g & \text{on }\partial \Omega .
\end{array}
\right. 
\end{equation*}
converge uniformly in $\bar{\Omega}$ to \emph{different} continuous limits $%
u_{0}$ and $u_{1}$.
\end{remark}

In \cite{HL}, Theorem \ref{FKP0} was extended to elliptic divergence form
operators with a singular drift:

\begin{theorem}[Theorem 2.17 in~\TeXButton{cite: HL}{\cite{HL}}]
\label{HL217} If $\mathcal{L}_{D,0}=\limfunc{div}\mathbf{A}_{0}\nabla +%
\mathbf{b}_{0}\cdot \nabla $ and $\mathcal{L}_{D,1}=\limfunc{div}\mathbf{A}%
_{1}\nabla +\mathbf{b}_{1}\cdot \nabla $ where $\mathbf{A}_{0}$ and $\mathbf{%
A}_{1}$ satisfy (\ref{CarlFKP}), and $\mathbf{b}_{i}=\left( b_{j}^{i}\right)
_{j=1}^{n+1}$, $i=0,1$ satisfy 
\begin{equation}
\delta \left( X\right) \sup_{Y\in Q_{\frac{\delta \left( X\right) }{2\sqrt{n}%
}}}\left| \mathbf{b}_{1}\left( Y\right) -\mathbf{b}_{0}\left( Y\right)
\right| ^{2}dX\in \frak{C}\left( d\sigma ,\Omega \right) ,  \label{CarlHL}
\end{equation}
then $d\omega _{D,0}\in A_{\infty }\Rightarrow d\omega _{D,1}\in A_{\infty }$%
.
\end{theorem}

The results in Theorems \ref{FKP0} and \ref{HL217} concern perturbation of
elliptic operators. They provide solvability for the $L^{q}$ Dirichlet
problem (for some $q>1$) for an operator $\mathcal{L}_{1}$ given that there
exists an operator $\mathcal{L}_{0}$ for which the $L^{p}$ Dirichlet problem
is solvable for some $p>1\,$\emph{and} the disagreement of their
coefficients satisfy the Carleson measure conditions (\ref{CarlFKP}) and (%
\ref{CarlHL}). In \cite{KP}, the authors answer a different question: What
are sufficient conditions on the coefficients $\mathbf{A}$ and $\mathbf{b}$
so that a given operator $\mathcal{L}_{D}=\limfunc{div}\mathbf{A}\nabla +%
\mathbf{b}$ has unique solutions for the $L^{p}$-Dirichlet problem for some $%
p>1$? See also \cite{KKPT}.

\begin{theorem}
\label{KP218}\cite{KP}Let $\mathcal{L}_{D}=\mathop{\rm div}\mathbf{A}\nabla +%
\mathbf{b}\cdot \nabla $, where $\mathbf{A}$ satisfies 
\begin{equation}
\delta \left( X\right) \sup_{\left| Y-X\right| \le \frac{\delta \left(
X\right) }{2}}\left| \nabla \mathbf{A}\left( Y\right) \right| ^{2}dX\in 
\frak{C}\left( d\sigma ,\Omega \right) ,  \label{CarlKP}
\end{equation}
and $\mathbf{b}$ satisfies (\ref{CarlHL}). Then $d\omega _{\mathcal{L}%
_{D}}\in A_{\infty }$.
\end{theorem}

\section{Statement of the results}

One of the main results in this work is an analog of Theorem \ref{HL217} for
nondivergence form operators.

\begin{theorem}
\label{CRth2}Let $\mathcal{L}_{N,0}=\mathbf{A}_{0}\cdot \nabla ^{2}+\mathbf{b%
}_{0}\cdot \nabla $ and $\mathcal{L}_{N,1}=\mathbf{A}_{1}\nabla ^{2}+\mathbf{%
b}_{1}\cdot \nabla $ where $\mathbf{A}_{\ell }=\left( A_{ij}^{\ell }\right)
_{i,j=1}^{n+1}$ and $\mathbf{b}_{\ell }=\left( b_{j}^{\ell }\right)
_{j=1}^{n+1}$ $\ell =0,1$ are bounded, measurable coefficients and $\mathbf{A%
}_{\ell }$ satisfy the ellipticity condition (\ref{ellipticity}) for $\ell
=0,1$. Suppose that $\mathcal{CD}$ holds for $\mathcal{L}_{N,\ell }$, $\ell
=0,1$ and that 
\begin{equation}
\left\{ \sup_{Y\in Q_{\frac{\delta \left( X\right) }{2\sqrt{n}}}\left(
X\right) }\frac{\left| \mathbf{A}_{1}\left( Y\right) -\mathbf{A}_{0}\left(
Y\right) \right| ^{2}+\delta ^{2}\left( X\right) \left| \mathbf{b}_{1}\left(
Y\right) -\mathbf{b}_{0}\left( Y\right) \right| ^{2}}{\delta \left( X\right) 
}\right\} dX\in \frak{C}\left( d\sigma ,\Omega \right) ,  \label{CarlCR}
\end{equation}
where $\sigma $\thinspace is a doubling measure on $\partial \Omega $. Then $%
d\omega _{N,0}\in A_{\infty }\Rightarrow d\omega _{N,1}\in A_{\infty }$.
\end{theorem}

As an application of this result, Theorem \ref{FKP0} (nondivergence case)
and an averaging of the coefficients argument we obtain an analog to Theorem 
\ref{KP218}. Moreover, this averaging argument (Lemma \ref{grcarl}) and
Theorem \ref{HL217} yield an extension of Theorem \ref{KP218} to the case
when condition (\ref{CarlKP}) is replaced by the weaker assumption (\ref
{CarlFKP}). To make the statements concise, we introduce the following
definition.

\begin{definition}
\label{oscillacions}For $r>0$, the $r$-oscillation of a measurable function $%
f\left( X\right) $ (scalar or vector-valued) at a point $X$ , denoted $%
\limfunc{osc}\nolimits_{r}f\left( X\right) $, is given by 
\begin{equation*}
\limfunc{osc}\nolimits_{r}f\left( X\right) =\sup_{Y,W\in Q_{r}\left(
X\right) }\left| f\left( W\right) -f\left( Z\right) \right| .
\end{equation*}
\end{definition}

\begin{theorem}
\label{CRth1}Let $\mathcal{L}_{D}=\limfunc{div}\mathbf{A}\nabla +\mathbf{b}%
\cdot \nabla $ and $\mathcal{L}_{N}=\mathbf{A}\cdot \nabla ^{2}+\mathbf{b}%
\cdot \nabla $ be uniformly elliptic operators in divergence form and
nondivergence form, respectively, with bounded measurable coefficient matrix 
$\mathbf{A}$ and drift vector $\mathbf{b}$. In the nondivergence case we
assume that $\mathcal{CD}$ holds for $\mathcal{L}_{N}$. Suppose that the
coefficients $\mathbf{A}$, $\mathbf{b}$ satisfy 
\begin{equation}
\frac{\left( \func{osc}_{\frac{\delta \left( X\right) }{2\sqrt{n}}}\mathbf{A}%
\left( X\right) \right) ^{2}+\delta ^{2}\left( X\right) \left( \func{osc}_{%
\frac{\delta \left( X\right) }{2\sqrt{n}}}\mathbf{b}\left( X\right) \right)
^{2}}{\delta \left( X\right) }dX\in \frak{C}\left( d\sigma ,\Omega \right) ,
\label{CarlCRa}
\end{equation}
Then $d\omega _{\mathcal{L}_{D}}\in A_{\infty }$ and $d\omega _{\mathcal{L}%
_{N}}\in A_{\infty }$.
\end{theorem}

\subsection{The theorems are sharp}

In \cite{FKP} it was shown that Theorem \ref{FKP0} is sharp for divergence
form equations in two fundamental ways (see Theorems 4.11 and 4.2 in \cite
{FKP}). The examples provided in that work were constructed using
Beurling-Ahlfors quasiconformal mappings on the half plane $\Bbb{R}%
_{+}^{2}\, $\cite{BA}. Quasi-conformal mappings preserve the divergence
form-structure of an elliptic operator $\mathcal{L}_{D}$, but when composed
with a nondivergence form operator $\mathcal{L}_{N}$ the transformed
operator has first order drift terms.

The more recent work \cite{KP} for divergence form operators (Theorem \ref
{KP218}) can be applied to obtain regularity of elliptic equation in \emph{%
nondivergence} form in the special case that the coefficient matrix
satisfies (\ref{CarlKP}). The approximation technique to be introduced in
the next section, together with Theorem \ref{CRth2}, show that Theorem \ref
{KP218} can be extended to operators in divergence and nondivergence form
with coefficients satisfying the weaker condition (\ref{CarlCRa}). At the
same time, this opens the door to extend the scope of the examples provided
in \ref{FKP0} to this wider class of operators. The following theorem is a
nondivergence analog to Theorem 4.11 in \cite{FKP}.

\begin{theorem}
\label{sharp}Given any nonnegative function $\alpha \left( X\right) $ in $%
\Bbb{R}_{+}^{2}=\left\{ \left( x,t\right) :t>0\right\} $ such that $\alpha
\left( X\right) $ satisfies the doubling condition: $\alpha \left( X\right)
\le C\alpha \left( X_{0}\right) $ for all $X=\left( x,t\right)
,\,X_{0}=\left( x_{0},t_{0}\right) :\,\left| X-X_{0}\right| <\frac{t_{0}}{2}$
and such that 
\begin{equation*}
\sup_{Y\in Q_{\frac{\delta \left( X\right) }{2\sqrt{n}}}\left( X\right) }%
\frac{\alpha \left( Y\right) ^{2}}{\delta \left( X\right) }dX\notin \frak{C}%
\left( d\sigma ,\left[ 0,1\right] ^{2}\right) ,
\end{equation*}
where $d\sigma $ is the Euclidean measure in $\partial \left( \left[
0,1\right] ^{2}\right) $ and $\delta \left( \left( x,t\right) \right) =t$.
There exists a coefficients matrix $\mathbf{A}$ such that, \newline
\emph{(1)}\quad the function $a\left( X\right) =\sup_{Y\in Q_{\frac{\delta
\left( X\right) }{2\sqrt{n}}}\left( X\right) }\left| \mathbf{A}\left(
Y\right) -\mathbf{I}\right| $, satisfies that for all $I\subset \Bbb{R}$, 
\begin{equation*}
\frac{1}{\left| I\right| }\int \int_{T\left( I\right) }a^{2}\left(
x,y\right) dx\,\frac{dy}{y}\le C\left[ \frac{1}{\left| I\right| }\int
\int_{T\left( 2I\right) }\alpha ^{2}\left( x,y\right) dx\,\frac{dy}{y}%
+1\right] ;
\end{equation*}
\newline
\emph{(2)}\quad the function $\tilde{a}\left( X\right) =\limfunc{osc}%
\nolimits_{Q_{\frac{\delta \left( X\right) }{2\sqrt{n}}}}\mathbf{A}\left(
X\right) $, satisfies that for all $I\subset \Bbb{R}$, 
\begin{equation*}
\frac{1}{\left| I\right| }\int \int_{T\left( I\right) }\tilde{a}^{2}\left(
x,y\right) dx\,\frac{dy}{y}\le C\left[ \frac{1}{\left| I\right| }\int
\int_{T\left( 2I\right) }\alpha ^{2}\left( x,y\right) dx\,\frac{dy}{y}%
+1\right] ;\qquad \text{and}
\end{equation*}
\newline
\emph{(3)}\quad if $\mathcal{L}_{N}=\mathbf{A}\cdot \nabla ^{2}$ on $\Bbb{R}%
_{+}^{2}$ , the elliptic measure $d\omega _{\mathcal{L}_{N}}$ is \emph{not}
in $A_{\infty }\left( dx,\left[ 0,1\right] \right) $ .\newline
\end{theorem}

The above theorem shows that the Carleson measure condition (\ref{CarlFKP})
in Theorem \ref{FKP0} is sharp also in the nondivergence case. In
particular, it shows that the main result in \cite{CR} is sharp. The proof
is a simple application of Theorem 4.11 in \cite{FKP} and the approximation
argument given in the section below (Lemma \ref{carlesin}). See the proof of
Theorem \ref{CRth2} in Section \ref{proofs} for an example of the techniques
used.

\section{Preliminary results}

Given a weight $w$\thinspace in the Muckenphout class $A_{p}\left( \Omega
\right) $, we denote by $L^{q}\left( \Omega ,w\right) ,$ $1\le p\le q<\infty 
$, the space of measurable functions $f$ such that 
\begin{equation*}
\left\| f\right\| _{L^{q}\left( \Omega ,w\right) }=\left( \int_{\Omega
}\left| f\left( x\right) \right| ^{q}w\left( x\right) dx\right) <\infty .
\end{equation*}
And for a nonnegative integer $k$, we define the Sobolev space $%
W^{k,q}\left( \Omega ,w\right) $ as the space of functions $f$ in $%
L^{q}\left( \Omega ,w\right) \,$such that $f$ has weak derivatives up to
order $k$ in $L^{q}\left( \Omega ,w\right) $. Under the assumption $w\in
A_{p}$, the space $W^{k,q}\left( \Omega ,w\right) $ is a Banach space and it
is also given as the closure of $\mathcal{C}_{0}^{\infty }\left( \Omega
\right) $ (smooth functions of compact support in $\Omega $) under the norm 
\begin{equation*}
\left\| f\right\| _{W^{k,q}\left( \Omega ,w\right) }=\sum_{\ell
=0}^{k}\left\| \nabla ^{\ell }f\right\| _{L^{q}\left( \Omega ,w\right) },
\end{equation*}
see \cite{FKS}, \cite{Ki}. We recall now some definitions.

\begin{definition}
Given a locally integrable function $f$ in $\Omega \subset \Bbb{R}^{n}$, the 
$\func{BMO}$ modulus of continuity of $f$ , $\eta _{\Omega ,f}\left(
r\right) $,is given by 
\begin{equation}
\eta _{\Omega ,f}\left( r\right) =\sup_{x\in \Omega }\sup_{0<s\le r}\frac{1}{%
\left| B_{s}\left( x\right) \bigcap \Omega \right| }\int_{B_{s}\left(
x\right) \bigcap \Omega }\left| f\left( y\right) -f_{s}\left( x\right)
\right| dy,  \label{bmo}
\end{equation}
where 
\begin{equation*}
f_{s}\left( x\right) =\frac{1}{\left| B_{s}\left( x\right) \bigcap \Omega
\right| }\int_{B_{s}\left( x\right) \bigcap \Omega }f\left( z\right) \,dz.
\end{equation*}
\end{definition}

The space $\limfunc{BMO}\left( \Omega \right) $ is given by 
\begin{equation*}
\limfunc{BMO}\left( \Omega \right) =\left\{ f\in L_{\limfunc{loc}}^{1}\left(
\Omega \right) :\left\| \eta _{\Omega ,f}\right\| _{L^{\infty }\left( \Bbb{R}%
_{+}\right) }<\infty \right\} .
\end{equation*}
For $\varrho \ge 0$, we let $\limfunc{BMO}\nolimits_{\varrho }\left( \Omega
\right) $ be given by 
\begin{equation*}
\limfunc{BMO}\nolimits_{\varrho }\left( \Omega \right) =\left\{ f\in 
\limfunc{BMO}\left( \Omega \right) :\liminf_{r\rightarrow 0^{+}}\eta
_{\Omega ,f}\left( r\right) \le \varrho \right\} .
\end{equation*}
It is easy to check that $\limfunc{BMO}\nolimits_{\varrho }\left( \Omega
\right) $ is a closed convex subset of $\limfunc{BMO}\left( \Omega \right) $
under the $\limfunc{BMO}$ norm $\left\| f\right\| _{\limfunc{BMO}\left(
\Omega \right) }=\left\| \eta _{\Omega ,f}\right\| _{L^{\infty }\left( \Bbb{R%
}_{+}\right) }$. When $\varrho =0,$ $\limfunc{BMO}\nolimits_{0}\left( \Omega
\right) $ is the space $\func{BMO}\left( \Omega \right) $ of functions of 
\emph{vanishing mean oscillation}. We say that a vector or matrix function
belongs to a space of scalar functions $\mathcal{X}$ if each component
belongs to that space $\mathcal{X}$. For example, we sat that the
coefficient matrix $\mathbf{A}\in \func{BMO}$ if $A_{ij}\in \func{BMO}$ for $%
1\le i,j\le n+1$. The following theorem establishes the solvability of the
continuous Dirichlet problem for a large class of elliptic operators in
nondivergence form. In particular, for such operators the harmonic measure
is well defined.

\begin{theorem}
\cite{CRt}\label{CRTth}Let $\mathcal{L}_{n}=a\cdot \nabla ^{2}$ be an
elliptic operator in nondivergence form in $\Omega \subset \Bbb{R}^{n+1}$
with ellipticity constants $\left( \lambda ,\Lambda \right) $. There exists
a constant $\varrho =\varrho \left( n,\lambda ,\Lambda \right) $, such that
if $a\in \limfunc{BMO}\nolimits_{\varrho }\left( \Omega \right) $ then for
every $g\in \mathcal{C}\left( \partial \Omega \right) $ there exists a
unique 
\begin{equation*}
u\in \mathcal{C}\left( \bar{\Omega}\right) \bigcap_{1\le p<\infty }W_{%
\limfunc{loc}}^{2,p}\left( \Omega \right)
\end{equation*}
such that $\mathcal{L}_{n}u=0$ in $\Omega $ and $u\equiv g$ on $\partial
\Omega $. Moreover, for any subdomain $\Omega ^{\prime }\subset \subset
\Omega $ and $1\le p<\infty $, there exists $C=C\left( n,\lambda ,\Lambda
,p,\Omega ,\limfunc{dist}\left( \Omega ^{\prime },\partial \Omega \right)
\right) $ such that the above solution satisfies 
\begin{equation*}
\left\| u\right\| _{W^{2,p}\left( \Omega ^{\prime }\right) }\le C\left\|
u\right\| _{L^{p}\left( \tilde{\Omega}^{\prime }\right) },\qquad \text{where 
}\tilde{\Omega}^{\prime }=\left\{ X\in \Omega :\delta \left( X\right) >\frac{%
1}{2}\limfunc{dist}\left( \Omega ^{\prime },\partial \Omega \right) \right\}
.
\end{equation*}
\end{theorem}

\begin{theorem}
{\label{strrong}}Let $w\in A_{p}$, $p\in [n,\infty )$ and $\Omega \subset 
\Bbb{R}^{n}$ be a Lipschitz domain. For any $0<\lambda \le \Lambda <\infty $
there exist a positive $\varrho =\varrho (n,p,\lambda ,\Lambda
,|[w]|_{A_{p}})$, such that if $\mathcal{L}_{N}\mathcal{=}\mathbf{A}\left(
X\right) \cdot \nabla _{X}^{2}$, $\mathbf{A}\in \limfunc{BMO}_{\rho }\left(
\Omega \right) $ then for any $f\in L^{p}(\Omega ,w)$, there exists a unique 
$u\in \mathcal{C}(\overline{\Omega })\bigcap W_{\limfunc{loc}}^{2,p}(\Omega
,w)$ such that $\mathcal{L}_{N}u=f$ in $\Omega $ and $u=0$ on $\partial
\Omega $. Moreover, if $\partial \Omega $ is of class $\mathcal{C}^{2}$,
then $u\in W_{\!0}^{1,p}(\Omega ,w)\bigcap W^{2,p}(\Omega ,w)$ and there
exists a positive $c=c(n,p,\lambda ,\Lambda ,|[w]|_{A_{p}},\eta _{\Omega ,%
\mathbf{A}},\partial \Omega )$, with $\eta _{\Omega ,\mathbf{A}}$ given by (%
\ref{bmo}), such that 
\begin{equation*}
\Vert u\Vert _{W^{2,p}(\Omega ,w)}\le c\,\Vert f\Vert _{L^{p}(\Omega ,w)}.
\end{equation*}
\end{theorem}

The following comparison principle is the main tool that allows us to treat
nondivergence form equations with singular drift. Since we assume that our
operator satisfies $\mathcal{CD}$, the result from \cite{BAR} originally
stated in a $\mathcal{C}^{2}$ domain and for continuous coefficients extends
to the more general case stated here. Before stating the theorem, we
introduce some notation.

If $\Omega \subset \Bbb{R}^{n+1}$ is a Lipschitz domain and $Q\in \partial
\Omega $, $r>0$, we define the \emph{boundary ball of radius }$r$\emph{\ at }%
$Q$ as 
\begin{equation*}
\triangle _{r}\left( Q\right) =\left\{ P\in \partial \Omega :\left|
P-Q\right| <r\right\} .
\end{equation*}
The \emph{Carleson region} associated to $\triangle _{r}\left( Q\right) \,$%
is 
\begin{equation*}
T_{r}\left( Q\right) =\left\{ X\in \Omega :\left| X-Q\right| <r\right\} .
\end{equation*}
The \emph{nontangential cone of aperture }$\alpha $\emph{\ and height }$r$%
\emph{\ at }$Q$ , $\alpha ,r>0$, is defined by 
\begin{equation*}
\Gamma _{\alpha ,r}\left( Q\right) =\left\{ X\in \Omega :\left| X-Q\right|
<\left( 1+\alpha \right) \delta \left( X\right) <\left( 1+\alpha \right)
r\right\} .
\end{equation*}

\begin{theorem}[Comparison Theorem for Solutions]
\cite{BAR}\label{comparison}Let $\Omega \subset \Bbb{R}^{n+1}$ be a
Lipschitz domain and $\mathcal{L}=\mathbf{A}\cdot \nabla ^{2}+\mathbf{b}%
\cdot \nabla $, be a uniformly elliptic operator such that $\mathbf{b}$
satisfies 
\begin{equation}
\delta \left( X\right) \sup_{Y\in Q_{\frac{\delta \left( X\right) }{2\sqrt{n}%
}}\left( X\right) }\left| \mathbf{b}\left( Y\right) \right| ^{2}dX\in \frak{C%
}\left( d\sigma ,\Omega \right)  \label{balone}
\end{equation}
and $\mathcal{L}$ satisfies $\mathcal{CD}$. There exists a constant $C$, $%
r_{0}>0$ depending only on $\mathcal{L}$ and $\Omega $ such that if $u$ and $%
v$ are two nonnegative solutions to $\mathcal{L}w=0$ in $T_{4r}\left(
Q\right) $, for some $Q\in \partial \Omega $, and such that $u\equiv v\equiv
0$ continuously on $\triangle _{2r}\left( Q\right) $, then 
\begin{equation*}
C^{-1}\frac{u\left( X\right) }{u\left( X_{r}\left( Q\right) \right) }\le 
\frac{v\left( X\right) }{v\left( X_{r}\left( Q\right) \right) }\le C\frac{%
u\left( X\right) }{u\left( X_{r}\left( Q\right) \right) },
\end{equation*}
for every $X\in T_{r}\left( Q\right) $, $0<r<r_{0}$. Here $X_{r}\left(
Q\right) \in \Omega $ with $\limfunc{dist}\left( X_{r}\left( Q\right)
\right) \approx \left| X-Q\right| \approx r$.
\end{theorem}

\begin{remark}
In \cite{BAR} the hypothesis on $\mathbf{b}$ is $\delta \left( X\right)
\left| \mathbf{b}\left( X\right) \right| \le \eta \left( \delta \left(
X\right) \right) $ where $\eta $ is an non-decreasing function such that $%
\eta \left( 0\right) =\lim_{s\rightarrow 0^{+}}\eta \left( s\right) =0$.
Nevertheless, the main tools used to obtain Theorem \ref{comparison} were%
\newline
\emph{(1)}\qquad $\mathcal{L}$ satisfies property $\mathcal{DP}$\newline
\emph{(2)}\qquad a maximum principle\newline
\emph{(3)}\qquad Harnack inequality.\newline
We assume \emph{(1)} holds, while $\emph{(2)}$ and \emph{(3)} follow as in 
\cite{BAR} once we notice that (\ref{balone}) implies that $\mathbf{b}$ is
locally bounded in $\Omega $. Therefore, Theorem \ref{comparison} also holds
under our assumptions.
\end{remark}

We list now some consequences of the above result that will be useful to us.

\begin{lemma}
\label{pro} Let $\Omega \subset \Bbb{R}^{n+1}$ be a Lipschitz domain and $%
\mathcal{L}=\mathbf{A}\cdot \nabla ^{2}+\mathbf{b}\cdot \nabla $, be a
uniformly elliptic operator such that $\mathbf{b}$ satisfies (\ref{balone})
and $\mathcal{L}$ satisfies $\mathcal{CD}$. Let $\triangle _{r}\left(
Q\right) $, $T_{r}\left( Q\right) $, $X_{r}\left( Q\right) \,$ and $r_{0}$
be as in Theorem \ref{comparison}. Let $\triangle =\triangle _{r}\left(
Q\right) $ for some $0<r\le r_{0}$ and $Q\in \partial \Omega $.\newline
\emph{(1)}\qquad $\omega ^{X_{r}\left( Q\right) }\left( \triangle \right)
\approx 1.$\newline
\emph{(2)}\qquad If $E\subset \triangle $ and $X\in \Omega \backslash
T_{2r}\left( Q\right) $, then 
\begin{equation*}
\omega ^{X_{r}\left( Q\right) }\left( E\right) \approx \frac{\omega
^{X}\left( E\right) }{\omega ^{X}\left( \triangle \right) }.
\end{equation*}
\newline
\emph{(3)}\qquad (Doubling property) $\omega ^{X}\left( \triangle _{r}\left(
Q\right) \right) \approx \omega ^{X}\left( \triangle _{2r}\left( Q\right)
\right) $ whenever $X\in \Omega \backslash T_{2r}\left( Q\right) $.
\end{lemma}

An important consequence of the properties listed in Lemma \ref{pro} is the
following analog of the ``main lemma'' in \cite{DJK}. Before stating the
result we need to introduce the concept of \emph{saw-tooth} region.

\begin{definition}
\label{serrucho}Given a Lipschitz domain $\Omega $ and $F\subset \partial
\Omega $ a closed set, a ``saw-tooth'' region $\Omega _{F}$ above $F$ in $%
\Omega $ of height $r>0$ is a Lipschitz subdomain of $\Omega $ with the
following properties:\newline
\emph{(1)}\qquad for some $0<\alpha <\beta $, $0<c_{1}<c_{2}$ and all $%
\alpha <\alpha ^{\prime }<\alpha ^{\prime \prime }<\beta ,$%
\begin{equation*}
\bigcup_{P\in F}\overline{\Gamma }_{\alpha ^{\prime },c_{1}r}\left( P\right)
\subset \Omega _{F}\subset \bigcup_{P\in F}\overline{\Gamma }_{\alpha
^{\prime \prime },c_{2}r}\left( P\right) ;
\end{equation*}
\newline
\emph{(2)}\qquad $\partial \Omega \bigcap \partial \Omega _{F}=F;$\newline
\emph{(3)}\qquad there exists $X_{0}\in \Omega _{F}$ (the center of $\Omega
_{F}$) such that $\limfunc{dist}\left( X_{0},\partial \Omega _{F}\right)
\approx r$,\newline
\emph{(4)}\qquad \label{also}for any $X\in \Omega $, $Q\in \partial \Omega $
such that $X\in \Gamma _{\alpha _{0},r_{0}}\left( Q\right) \bigcap \Omega
_{F}\neq \emptyset $, $\exists P\in F:Q_{\frac{\delta \left( X\right) }{2%
\sqrt{n}}}\left( X\right) \subset \Gamma _{\alpha _{0},r_{0}}\left( P\right)
,$\newline
\emph{(5)}\qquad $\Omega _{F}$\thinspace is a Lipschitz domain with
Lipschitz constant depending only on that of $\Omega $.
\end{definition}

With these provisos, now we state the following:

\begin{lemma}
\label{mmain}Let $\Omega $ be a Lipschitz domain and $\mathcal{L}=\mathbf{A}%
\cdot \nabla ^{2}+\mathbf{b}\cdot \nabla $, be a uniformly elliptic operator
that satisfies $\mathcal{CD}$. Let $F\subset \partial \Omega $ be a closed
set, and let $\Omega _{F}$ be a saw-tooth region above $F$ in $\Omega $. Let 
$\omega =\omega _{\mathcal{L},\Omega }$ and let $\nu =\omega _{\mathcal{L}%
,\Omega _{F}}^{X_{0}}$ , where $X_{0}$ is the center of $\Omega _{F}$. There
exists $\theta >0$ such that 
\begin{equation*}
\frac{\omega \left( E\right) }{\omega \left( \triangle \right) }\le \nu
\left( F\right) ^{\theta },\qquad \text{for }F\subset \triangle =\triangle
_{r}\left( Q\right) .
\end{equation*}
Here $\theta $ depends on the Lipschitz character of $\Omega $, but not on $%
E $ or $\triangle $.
\end{lemma}

The following result will allow us to relax the hypothesis (\ref{CarlKP}) in
Theorem \ref{KP218}. We adopt the following notation 
\begin{eqnarray*}
W &=&\left( w,r\right) ,\qquad \text{where}\quad w=\left( w_{1},\dots
,w_{n}\right) ,\qquad r=w_{n+1} \\
Y &=&\left( y,s\right) ,\qquad \text{where}\quad y=\left( y_{1},\dots
,y_{n}\right) ,\qquad s=y_{n+1} \\
X &=&\left( x,t\right) ,\qquad \text{where}\quad x=\left( x_{1},\dots
,x_{n}\right) ,\qquad t=x_{n+1}.
\end{eqnarray*}
Given a measurable function $f\left( X\right) $ and $\gamma >0$, we recall
that the oscillation of $f$ in the cube $Q_{\gamma }\left( X\right) $ is
given by 
\begin{equation*}
\limfunc{osc}\nolimits_{\gamma }f\left( X\right) =\sup_{Y,W\in Q_{\gamma
}\left( X\right) }\left| f\left( W\right) -f\left( Z\right) \right| .
\end{equation*}

\begin{lemma}
\label{grcarl}Let $f\left( X\right) $ be a bounded measurable function in $%
\Bbb{R}_{+}^{n+1}$. Let $0<d_{0}<1$ and let $a,b,c$ be positive constants
such that $a+b+c\le d_{0}$. Define 
\begin{eqnarray}
\tilde{f}\left( Y\right) &=&\frac{1}{\left| Q_{a}\left( Y\right) \right| }%
\int_{Q_{a}\left( Y\right) }f\left( W\right) dW  \label{ftilde} \\
f^{*}\left( X\right) &=&\frac{1}{\left| Q_{b}\left( X\right) \right| }%
\int_{Q_{b}\left( X\right) }\tilde{f}\left( Y\right) dY.  \notag
\end{eqnarray}
Then, if $\mathcal{E}\left( X\right) =\limfunc{osc}\nolimits_{d_{0}}f\left(
X\right) $, $\mathcal{E}^{*}\left( X\right) =\sup_{Y\in Q_{c}\left( X\right)
}\left| \nabla f^{*}\left( Y\right) \right| $ and $\mathcal{\tilde{E}}\left(
X\right) =\sup_{Y\in Q_{c}\left( X\right) }\left| f^{*}\left( Y\right)
-f\left( Y\right) \right| $, they satisfy 
\begin{equation*}
\mathcal{E}^{*}\left( X\right) \le C_{a,b,n}\frac{\mathcal{E}\left( X\right) 
}{t}\qquad \text{and}\qquad \mathcal{\tilde{E}}\left( X\right) \le \mathcal{E%
}\left( X\right) ,
\end{equation*}
respectively. In particular, if $\mu $ is a doubling measure on $\partial 
\Bbb{R}_{+}^{n+1}\,$and 
\begin{equation*}
\frac{\mathcal{E}\left( X\right) ^{2}}{t}dX\in \frak{C}\left( d\mu ,\Bbb{R}%
_{+}^{n+1}\right) \text{, with norm }\left\| \frac{\mathcal{E}\left(
X\right) ^{2}}{t}\right\| _{\frak{C}\left( d\mu ,\Bbb{R}_{+}^{n+1}\right) }=%
\mathcal{C}_{0}.
\end{equation*}
then $t\mathcal{E}^{*}\left( X\right) ^{2}dX\in \frak{C}\left( d\mu ,\Bbb{R}%
_{+}^{n+1}\right) $ (with norm $\mathcal{C}^{*}=\mathcal{C}^{*}\left( 
\mathcal{C}_{0},n,a,b\right) $) and $t^{-1}\mathcal{\tilde{E}}\left(
X\right) ^{2}dX\in \frak{C}\left( d\mu ,\Bbb{R}_{+}^{n+1}\right) $.
\end{lemma}

$.$%
\proof%
For $\gamma >1$, let $\mathcal{T}_{i}^{\gamma }\left( Y\right) $ denote the
intersection of $Q_{\gamma }\left( Y\right) $ and the hyperplane normal to
the $i^{th}$ coordinate axis, i.e. 
\begin{equation*}
\mathcal{T}_{i}^{\gamma }\left( Y\right) =\left\{ Z\in Q_{\gamma }\left(
Y\right) :Z_{i}=Y_{i}\right\} ,\qquad i=1,\dots ,n+1.
\end{equation*}
Let $\vec{e}_{i}=\left( \delta _{ij}\right) _{j=1}^{n+1}$, $i=1,\dots ,n+1$,
and let $\hat{w}^{i}=\left( w_{1},\dots ,w_{i-1},0,w_{i+1},\dots
,w_{n},r\right) $, $d\hat{w}^{i}=dw_{1}\dots dw_{i-1}dw_{i+1}\dots dw_{n+1}$
with similar definitions for the variable $Y$. Then for $i=1,\dots ,n,$%
\begin{eqnarray}
\frac{\partial \tilde{f}\left( Y\right) }{\partial y_{i}} &=&\lim_{h%
\rightarrow 0}\frac{1}{h}\left\{ \frac{1}{\left| Q_{a}\left( Y+h\vec{e}%
_{i}\right) \right| }\int_{Q_{a}\left( Y+h\vec{e}_{i}\right) }fdW-\frac{1}{%
\left| Q_{a}\left( Y\right) \right| }\int_{Q_{a}\left( Y\right) }fdW\right\}
\notag \\
&=&\frac{1}{\left| Q_{a}\left( Y\right) \right| }\lim_{h\rightarrow 0}\frac{1%
}{h}\left\{ \int_{y_{i}+\frac{as}{2}}^{y_{i}+\frac{as}{2}+h}\int_{\mathcal{T}%
_{i}^{a}\left( Y\right) }fdW-\int_{y_{i}-\frac{as}{2}}^{y_{i}-\frac{as}{2}%
+h}\int_{\mathcal{T}_{i}^{a}\left( Y\right) }fdW\right\}  \notag \\
&=&\frac{1}{at}\frac{1}{\left| \mathcal{T}_{i}^{a}\left( Y\right) \right| }%
\left\{ \int_{\mathcal{T}_{i}^{a}\left( Y\right) }f\left( \hat{w}^{i}+\left(
y_{i}+\frac{as}{2}\right) \vec{e}_{i}\right) d\hat{w}^{i}\right.  \label{ftx}
\\
&&\left. -\int_{\mathcal{T}_{i}^{a}\left( Y\right) }f\left( \hat{w}%
^{i}+\left( y_{i}-\frac{as}{2}\right) \vec{e}_{i}\right) d\hat{w}%
^{i}\right\} .  \notag
\end{eqnarray}
And 
\begin{eqnarray}
\frac{\partial \tilde{f}\left( Y\right) }{\partial s} &=&\lim_{h\rightarrow
0}\frac{1}{h}\left\{ \frac{1}{\left| Q_{a}\left( Y+h\vec{e}_{n+1}\right)
\right| }\int_{Q_{a}\left( Y+h\vec{e}_{n+1}\right) }fdW-\frac{1}{\left|
Q_{a}\left( Y\right) \right| }\int_{Q_{a}\left( Y\right) }fdW\right\}  \notag
\\
&=&\lim_{h\rightarrow 0}\frac{1}{h}\left\{ \frac{1}{\left| Q_{a}\left( Y+h%
\vec{e}_{n+1}\right) \right| }\int_{Q_{a}\left( Y+h\vec{e}_{n+1}\right)
\backslash Q_{a}\left( Y\right) }fdW\right.  \notag \\
&&\left. -\frac{1}{\left| Q_{a}\left( Y\right) \right| }\int_{Q_{a}\left(
Y\right) \backslash Q_{a}\left( Y+h\vec{e}_{n+1}\right) }fdW\right\}  \notag
\\
&=&\frac{1}{as}\frac{1}{\left| \mathcal{T}_{n+1}^{a}\left( Y\right) \right| }%
\left\{ \int_{\mathcal{T}_{n+1}^{a}\left( Y\right) }\left[ f\left( w,s\left(
1+\frac{a}{2}\right) \right) -f\left( w,s\left( 1-\frac{a}{2}\right) \right)
\right] dw\right\}  \notag \\
&&+\frac{1}{as}\sum_{i=1}^{n}\frac{1}{\left| \mathcal{T}_{i}^{a}\left(
Y\right) \right| }\int_{\mathcal{T}_{i}^{a}\left( Y\right) }f\left( \hat{w}%
^{i}+\left( y_{i}+\frac{as}{2}\right) \vec{e}_{i}\right) d\hat{w}^{i}ds 
\notag \\
&&+\frac{1}{as}\sum_{i=1}^{n}\frac{1}{\left| \mathcal{T}_{i}^{a}\left(
Y\right) \right| }\int_{\mathcal{T}_{i}^{a}\left( Y\right) }f\left( \hat{w}%
^{i}+\left( y_{i}-\frac{as}{2}\right) \vec{e}_{i}\right) d\hat{w}^{i}ds 
\notag \\
&=&\frac{1}{as}\left\{ \tilde{F}_{1}\left( Y\right) +\tilde{F}_{2}\left(
Y\right) +\tilde{F}_{3}\left( Y\right) \right\} .  \label{ftt}
\end{eqnarray}
Now, let $Q$ be the unit cube $Q=\left\{ Z=\left( z_{1},\dots
,z_{n},r\right) :\left| z_{i}\right| <\frac{1}{2}\right\} $; substituting $Z=%
\frac{Y-X}{bt}$ we have 
\begin{equation*}
f^{*}\left( X\right) =\int_{Q}\tilde{f}\left( X+Zbt\right) dZ.
\end{equation*}
Hence, for $i=1,\dots ,n$%
\begin{eqnarray}
&&\frac{\partial f^{*}\left( X\right) }{\partial x_{i}}=\int_{Q}\frac{%
\partial }{\partial x_{i}}\tilde{f}\left( X+Zbt\right) dZ  \notag \\
&=&\frac{1}{\left| Q_{b}\right| }\int_{Q_{b}}\frac{\partial \tilde{f}\left(
Y\right) }{\partial x_{i}}dY  \notag \\
&=&\frac{1}{bt}\frac{1}{\left| \mathcal{T}_{i}^{b}\left( X\right) \right| }%
\left\{ \int_{\mathcal{T}_{i}^{b}\left( X\right) }\tilde{f}\left( \hat{y}%
^{i}+\left( x_{i}+\frac{bt}{2}\right) \vec{e}_{i}\right) d\hat{y}%
^{i}ds\right.  \label{Os} \\
&&\left. -\int_{\mathcal{T}_{i}^{b}\left( X\right) }\tilde{f}\left( \hat{y}%
^{i}+\left( x_{i}-\frac{bt}{2}\right) \vec{e}_{i}\right) d\hat{y}%
^{i}ds\right\} .  \notag
\end{eqnarray}
And 
\begin{eqnarray}
\frac{\partial f^{*}\left( X\right) }{\partial t} &=&\int_{Q}\frac{\partial 
}{\partial t}\tilde{f}\left( X+Zbt\right) dZ  \notag \\
&=&\sum_{i=1}^{n}\frac{1}{\left| Q_{b}\left( X\right) \right| }%
\int_{Q_{b}\left( X\right) }\frac{\left( y_{i}-x_{i}\right) }{t}\frac{%
\partial \tilde{f}\left( Y\right) }{\partial x_{i}}dY  \notag \\
&&+\frac{1}{\left| Q_{b}\left( X\right) \right| }\int_{Q_{b}\left( X\right)
}\left( 1+\frac{\left( s-t\right) }{t}\right) \frac{\partial \tilde{f}\left(
Y\right) }{\partial t}dY  \notag \\
&=&\sum_{i=1}^{n}I_{i}+I_{n+1}.  \label{Is}
\end{eqnarray}
Let $\sigma =\frac{b}{\ln \left( \frac{2+b}{2-b}\right) }$, then 
\begin{eqnarray}
&&\frac{1}{\left| Q_{b}\left( X\right) \right| }\int_{Q_{b}\left( X\right)
}\left( \frac{\left( s-t\right) }{t}+\left( 1-\sigma \right) \right) \frac{1%
}{s}dY  \notag \\
&=&\frac{1}{bt}\frac{1}{\left| \mathcal{T}_{n+1}^{b}\left( X\right) \right| }%
\int_{t-\frac{bt}{2}}^{t+\frac{bt}{2}}\int_{\mathcal{T}_{n+1}^{b}\left(
X\right) }\left( \frac{1}{t}-\frac{\sigma }{s}\right) dyds  \notag \\
&=&\frac{1}{bt}\int_{t-\frac{bt}{2}}^{t+\frac{bt}{2}}\left( \frac{1}{t}-%
\frac{\sigma }{s}\right) dyds  \notag \\
&=&\frac{1}{bt}\left( b-\sigma \ln \left( \frac{2+b}{2-b}\right) \right) =0.
\label{zerop}
\end{eqnarray}
Substituting $1+\frac{\left( s-t\right) }{t}=\sigma +\left( \frac{\left(
s-t\right) }{t}+\left( 1-\sigma \right) \right) $ on the right of (\ref{Is}%
), we have 
\begin{eqnarray}
I_{n+1} &=&\frac{\sigma }{\left| Q_{b}\left( X\right) \right| }%
\int_{Q_{b}\left( X\right) }\frac{\partial \tilde{f}\left( Y\right) }{%
\partial t}dY+\frac{1}{\left| Q_{b}\left( X\right) \right| }%
\int_{Q_{b}\left( X\right) }\left( \frac{\left( s-t\right) }{t}+\left(
1-\sigma \right) \right) \frac{\partial \tilde{f}\left( Y\right) }{\partial t%
}dY  \notag \\
&=&\frac{\sigma }{bt}\frac{1}{\left| \mathcal{T}_{n+1}^{b}\left( X\right)
\right| }\int_{\mathcal{T}_{n+1}^{b}\left( X\right) }\left[ \tilde{f}\left(
y,t+\frac{bt}{2}\right) -\tilde{f}\left( y,t-\frac{bt}{2}\right) \right] dy+%
\tilde{I}_{n+1}.  \label{In1}
\end{eqnarray}
Using (\ref{ftt}) in the second integral on the right, we have 
\begin{eqnarray}
\tilde{I}_{n+1} &=&\frac{1}{\left| Q_{b}\left( X\right) \right| }%
\int_{Q_{b}\left( X\right) }\left( \frac{\left( s-t\right) }{t}+\left(
1-\sigma \right) \right) \frac{1}{as}\tilde{F}_{1}dY  \label{Inp1} \\
&&+\frac{1}{\left| Q_{b}\left( X\right) \right| }\int_{Q_{b}\left( X\right)
}\left( \frac{\left( s-t\right) }{t}+\left( 1-\sigma \right) \right) \frac{1%
}{as}\left( \tilde{F}_{2}+\tilde{F}_{3}\right) dY  \notag
\end{eqnarray}
For $Y\in Q_{b}\left( X\right) $ and $W\in Q_{a}\left( Y\right) $, we have $%
\left| f\left( W\right) -f\left( X\right) \right| \le \limfunc{osc}%
\nolimits_{\left( a+b\right) }f\left( X\right) $. Then 
\begin{eqnarray}
\left| \tilde{F}_{1}\left( Y\right) \right| &\le &\frac{1}{\left| \mathcal{T}%
_{n+1}^{a}\left( Y\right) \right| }\int_{\mathcal{T}_{n+1}^{a}\left(
Y\right) }\left| f\left( w,s\left( 1+\frac{1}{2a}\right) \right) -f\left(
X\right) \right| dw  \notag \\
&&+\frac{1}{\left| \mathcal{T}_{n+1}^{a}\left( Y\right) \right| }\int_{%
\mathcal{T}_{n+1}^{a}\left( Y\right) }\left| f\left( w,s\left( 1-\frac{1}{2a}%
\right) \right) -f\left( X\right) \right| dw  \notag \\
&\le &2\left( \limfunc{osc}\nolimits_{a+b}f\right) \left( X\right) .
\label{f1t}
\end{eqnarray}
Also, 
\begin{eqnarray}
&&\left| \tilde{F}_{2}\left( Y\right) +\tilde{F}_{3}\left( Y\right)
-2nf\left( X\right) \right|  \notag \\
&=&\left| \sum_{i=1}^{n}\frac{1}{\left| \mathcal{T}_{i}^{a}\left( Y\right)
\right| }\left\{ \int_{\mathcal{T}_{i}^{a}\left( Y\right) }\left[ f\left( 
\hat{w}^{i}+\left( y_{i}+\frac{as}{2}\right) \vec{e}_{i}\right) -f\left(
X\right) \right] d\hat{w}^{i}ds\right\} \right.  \notag \\
&&\left. \sum_{i=1}^{n}\frac{1}{\left| \mathcal{T}_{i}^{a}\left( Y\right)
\right| }\left\{ \int_{\mathcal{T}_{i}^{a}\left( Y\right) }\left[ f\left( 
\hat{w}^{i}+\left( y_{i}-\frac{as}{2}\right) \vec{e}_{i}\right) -f\left(
X\right) \right] d\hat{w}^{i}ds\right\} \right|  \notag \\
&\le &2n\limfunc{osc}\nolimits_{\left( a+b\right) }f\left( X\right) .
\label{f2t}
\end{eqnarray}
Using the cancelation property (\ref{zerop}), (\ref{f1t}) and (\ref{f2t}) on
the right of (\ref{Inp1}), we obtain 
\begin{eqnarray}
\left| \tilde{I}_{n+1}\right| &\le &\frac{1}{at}\frac{1}{\left| Q_{b}\left(
X\right) \right| }\int_{Q_{b}\left( X\right) }\left| 1-\frac{t\mu }{s}%
\right| \left| \tilde{F}_{1}\right| dY  \notag \\
&&+\frac{1}{at}\frac{1}{\left| Q_{b}\left( X\right) \right| }%
\int_{Q_{b}\left( X\right) }\left| 1-\frac{t\mu }{s}\right| \left| \tilde{F}%
_{2}+\tilde{F}_{3}-2nf\left( X\right) \right| dY  \notag \\
&\le &\frac{4\left( n+1\right) }{at}\limfunc{osc}\nolimits_{\left(
a+b\right) }f\left( X\right) .  \label{Itnp1}
\end{eqnarray}
Where we also used that for $Y\in Q_{b}\left( X\right) $ and $b<1$, $\left|
1-\frac{t\mu }{s}\right| \le 2$. Substituting on the second term on the
right of (\ref{In1}), and handling the first term on the right of (\ref{In1}%
) in a similar manner, we obtain 
\begin{equation}
\left| I_{n+1}\right| \le \frac{\sigma }{bt}\limfunc{osc}\nolimits_{b}f%
\left( X\right) +\frac{4\left( n+1\right) }{at}\limfunc{osc}%
\nolimits_{\left( a+b\right) }f\left( X\right) \le \frac{\mu
b^{-1}+4a^{-1}\left( n+1\right) }{t}\limfunc{osc}\nolimits_{\left(
a+b\right) }f\left( X\right) .  \label{IInp1}
\end{equation}
The estimate for $I_{i}$, $i=1,\dots ,n$ in (\ref{Is}) is obtained also
proceeding in this way. Note that from (\ref{ftx}) we have 
\begin{equation*}
\left| \frac{\partial \tilde{f}\left( Y\right) }{\partial y_{i}}\right| \le 
\frac{1}{at}\limfunc{osc}\nolimits_{\left( a+b\right) }f\left( X\right) ,
\end{equation*}
hence 
\begin{equation*}
\left| I_{i}\right| \le \frac{1}{at}\limfunc{osc}\nolimits_{\left(
a+b\right) }f\left( X\right) \frac{1}{\left| Q_{b}\left( X\right) \right| }%
\int_{Q_{b}\left( X\right) }\left| \frac{y_{i}-x_{i}}{t}\right| dY\le \frac{1%
}{abt}\limfunc{osc}\nolimits_{\left( a+b\right) }f\left( X\right) .
\end{equation*}
Using these estimates and (\ref{IInp1}) on the right of (\ref{Is}) finally
yields 
\begin{equation*}
\left| \frac{\partial f^{*}\left( X\right) }{\partial t}\right| \le \frac{%
na^{-1}b^{-1}+\mu b^{-1}+4a^{-1}\left( n+1\right) }{t}\limfunc{osc}%
\nolimits_{\left( a+b\right) }f\left( X\right) =C_{a,b,n}\frac{\func{osc}%
_{\left( a+b\right) }f\left( X\right) }{t}.
\end{equation*}
On the other hand, from (\ref{Os}) and the definition of $\tilde{f}\left(
Y\right) $, (\ref{ftilde}), it follows that for $i=1,\dots ,n$, 
\begin{eqnarray*}
\left| \frac{\partial f^{*}\left( X\right) }{\partial x_{i}}\right| &\le &%
\frac{1}{bt}\frac{1}{\left| \mathcal{T}_{i}^{b}\left( X\right) \right| }%
\int_{\mathcal{T}_{i}^{b}\left( X\right) }\left| \tilde{f}\left( \hat{y}%
^{i}+\left( x_{i}+\frac{bt}{2}\right) \vec{e}_{i}\right) -\tilde{f}\left( 
\hat{y}^{i}+\left( x_{i}-\frac{bt}{2}\right) \vec{e}_{i}\right) \right| d%
\hat{y}^{i}ds \\
&\le &\frac{1}{bt}\limfunc{osc}\nolimits_{\left( a+b\right) }f\left(
X\right) .
\end{eqnarray*}
Hence, we have shown that 
\begin{equation*}
\left| \nabla f^{*}\left( X\right) \right| \le C_{a,b,n}\frac{\func{osc}%
_{\left( a+b\right) }f\left( X\right) }{t}.
\end{equation*}
And therefore, for $0<a,b,c<1$, satisfying $a+b+c\le d_{0}$ we have 
\begin{equation*}
\mathcal{E}^{*}\left( X\right) =\sup_{Y\in Q_{c}\left( X\right) }\left|
\nabla f^{*}\left( Y\right) \right| \le C_{a,b,n}\frac{\limfunc{osc}%
\nolimits_{\left( a+b+c\right) }f\left( X\right) }{t}\le C_{a,b,n}\frac{%
\mathcal{E}\left( X\right) }{t},
\end{equation*}
Then $t\mathcal{E}^{*}\left( X\right) ^{2}\le C_{a,b,n}^{2}\frac{\mathcal{E}%
^{2}\left( X\right) }{t}$ and if the right hand side satisfies $\frac{%
\mathcal{E}^{2}\left( X\right) }{t}\in \frak{C}\left( d\mu ,\Bbb{R}%
_{+}^{n+1}\right) $ with Carleson $\mathcal{C}_{0}$, it follows that $t%
\mathcal{E}^{*}\left( X\right) ^{2}\in \frak{C}\left( d\mu ,\Bbb{R}%
_{+}^{n+1}\right) $ with Carleson norm bounded by $C_{a,b,n}^{2}\mathcal{C}%
_{0}$. Finally, for $Y\in Q_{c}\left( X\right) $, 
\begin{eqnarray*}
\left| f^{*}\left( Y\right) -f\left( Y\right) \right| &=&\left| \frac{1}{%
\left| Q_{b}\left( Y\right) \right| }\int_{Q_{b}\left( Y\right) }\tilde{f}%
\left( Z\right) dZ-f\left( Y\right) \right| \\
&=&\left| \frac{1}{\left| Q_{b}\left( Y\right) \right| }\int_{Q_{b}\left(
Y\right) }\frac{1}{\left| Q_{a}\left( Z\right) \right| }\int_{Q_{a}\left(
Z\right) }f\left( W\right) dWdZ-f\left( Y\right) \right| \\
&\le &\frac{1}{\left| Q_{b}\left( Y\right) \right| }\int_{Q_{b}\left(
Y\right) }\frac{1}{\left| Q_{a}\left( Z\right) \right| }\int_{Q_{a}\left(
Z\right) }\left| f\left( W\right) -f\left( Y\right) \right| dWdZ \\
&\le &\limfunc{osc}\nolimits_{a+b+c}f\left( X\right) \le \mathcal{E}\left(
X\right) .
\end{eqnarray*}
Hence $\sup_{Y\in Q_{c}\left( X\right) }\left| f^{*}\left( Y\right) -f\left(
Y\right) \right| \le \mathcal{E}\left( X\right) $.%
\endproof%

\begin{lemma}
\label{carlesin}Let $\Omega $ be a Lipschitz domain in $\Bbb{R}^{n+1}$ , $%
\mu $ be a doubling measure on $\partial \Omega $, $\delta \left( X\right) =%
\limfunc{dist}\left( X,\partial \Omega \right) $, and suppose that $g$ is a
measurable function (scalar or vector valued) such that for some constant $%
0<\alpha <1$ it satisfies 
\begin{equation*}
d\nu =\delta \left( X\right) \limfunc{osc}\nolimits_{\alpha \delta \left(
X\right) }g\left( X\right) ^{2}dX\in \frak{C}\left( d\mu ,\Omega \right) .
\end{equation*}
The for every Lipschitz subdomain $\tilde{\Omega}\subset \Omega $, there
exists a doubling measure $\tilde{\mu}$ on $\partial \tilde{\Omega}$, with
doubling constant depending only on the doubling constant of $\mu $, such
that $\tilde{\mu}=\mu $ in $\partial \Omega \bigcap \partial \tilde{\Omega}$%
, and 
\begin{equation*}
d\tilde{\nu}=\tilde{\delta}\left( X\right) \limfunc{osc}\nolimits_{\alpha
\delta \left( X\right) }g\left( X\right) ^{2}dX\in \frak{C}\left( d\tilde{\mu%
},\tilde{\Omega}\right) .
\end{equation*}
where $\tilde{\delta}\left( X\right) =\limfunc{dist}\left( X,\partial \tilde{%
\Omega}\right) $. Moreover, the Carleson norm $\left\| d\tilde{\nu}\right\|
_{\frak{C}\left( d\tilde{\mu},\tilde{\Omega}\right) }$depends only on the
Carleson norm of $d\nu $ and the Lipschitz character of $\tilde{\Omega}$. If 
$\mu \,$is the Lebesgue measure on $\partial \Omega $ then $\tilde{\mu}$ can
be taken as the Lebesgue measure on $\partial \tilde{\Omega}$.
\end{lemma}

\proof%
In the case that $d\sigma $ is the Euclidean measure $d\sigma $ on $\partial
\Omega $, the result follows as in the proof of Lemma 3.1 in \cite{KP}. To
obtain the Lemma \ref{carlesin} for an arbitrary doubling measure $\sigma $,
given a Lipschitz subdomain $\tilde{\Omega}\subset \Omega $, we will
construct an appropriate extension of $\sigma $ from $\partial \Omega
\bigcap \partial \tilde{\Omega}$ to $\partial \tilde{\Omega}$. We may assume
that $\Omega =\Bbb{R}_{+}^{n+1}=\left\{ \left( x,t\right) :x\in \Bbb{R}%
^{n},\,t>0\right\} $, the general case follows by standard techniques. For $%
X=\left( x,t\right) \in \Bbb{R}_{+}^{n+1}$, let $\triangle _{X}=\left\{
\left( y,0\right) :\left| x-y\right| \le t\right\} $ and define $M\left(
X\right) =\frac{\sigma \left( \triangle _{X}\right) }{\sigma \left(
\triangle _{x}\right) }$, where $\sigma $ is the Euclidean measure in $\Bbb{R%
}^{n}=\partial \Omega $. Since $\sigma $ is doubling, $M\left( X\right) $ is
continuous in $\Bbb{R}_{+}^{n+1}$ and $M\rightarrow \frac{d\sigma }{d\sigma }
$ as $t\rightarrow 0$ in the weak$^{*}$ topology of measures. Now, for any
Borel set $E\subset \tilde{\Omega}$, we define 
\begin{equation*}
\tilde{\mu}\left( E\right) =\sigma \left( E\bigcap \partial \Omega \right)
+\int_{E\bigcap \partial \tilde{\Omega}}M\left( X\right) d\tilde{\sigma}%
\left( X\right) ,
\end{equation*}
where $d\tilde{\sigma}\left( X\right) $ denotes the Euclidean measure on $%
\partial \tilde{\Omega}$. Then obviously $\tilde{\mu}$ is an extension of $%
\sigma $ from $\partial \Omega \bigcap \partial \tilde{\Omega}$ to $\partial 
\tilde{\Omega}$. It remains to check that $\tilde{\mu}$ is a doubling
measure and that $d\tilde{\nu}\in \frak{C}\left( d\tilde{\mu},\tilde{\Omega}%
\right) $. Let $\triangle \subset \partial \tilde{\Omega}\mathcal{\ }$ be a
surface ball centered at $X_{0}=\left( x_{0},t_{0}\right) \in \partial 
\tilde{\Omega}$ and let $T\left( \triangle \right) \subset \tilde{\Omega}$
be the associated Carleson region. Following the ideas in \cite{KP}, we
consider two cases.

In case 1, we assume that $d_{0}=\mathop{\rm diam}T\left( \triangle \right) $
is smaller than $\frac{t_{0}}{10}$. If $X=\left( x,t\right) \in \Bbb{R}%
_{+}^{n+1}$ belongs to $T\left( \triangle \right) $, then $\left|
t_{0}-t\right| \le \frac{t_{0}}{10}$ and so $\delta \left( X_{0}\right)
=t_{0}\le \frac{10}{9}t=\frac{10}{9}\delta \left( X\right) $. From the
definition of $M\left( X\right) $, and the doubling property of $\sigma $,
it follows that for all $X\in B_{\frac{t_{0}}{10}}\left( X_{0}\right) $, $%
M\left( X\right) \approx M\left( X_{0}\right) $ with constants depending
only on $\sigma $. Since $\tilde{\Omega}$ is a Lipschitz domain, it follows
that 
\begin{equation*}
\tilde{\mu}\left( \triangle \right) =\int_{\triangle }M\left( X\right) d%
\tilde{\sigma}\left( X\right) \approx M\left( X_{0}\right) \int_{\triangle }d%
\tilde{\sigma}\left( X\right) \approx d_{0}^{n}M\left( X_{0}\right) \approx 
\frac{d_{0}^{n}}{t_{0}^{n}}\sigma \left( \triangle _{X_{0}}\right) .
\end{equation*}
Similarly, $\tilde{\mu}\left( 2\triangle \right) \approx \frac{2^{n}d_{0}^{n}%
}{t_{0}^{n}}\sigma \left( 2\triangle _{X_{0}}\right) $, which shows that $%
\tilde{\mu}\left( \triangle \right) \approx \tilde{\mu}\left( 2\triangle
\right) $. On the other hand, because $d\nu \left( X\right) \in \frak{C}%
\left( d\sigma ,\Omega \right) $, we have 
\begin{equation*}
\sup_{Y,Z\in B_{\alpha \tilde{\delta}\left( X\right) }\left( X\right)
}\left| g\left( Y\right) -g\left( Z\right) \right| \le \frac{C}{t}\left( 
\frac{\sigma \left( \triangle _{X}\right) }{t^{n}}\right) ^{\frac{1}{2}},
\end{equation*}
where $C$ depends only on $n$, $\alpha $, the doubling constant of $\sigma $
and $\left\| d\nu \right\| _{\frak{C}\left( d\sigma ,\Omega \right) }$.
Since $\tilde{\Omega}\subset \Omega $,we have $\tilde{\delta}\left( X\right)
\le \delta \left( X\right) =t$ for all $X\in \tilde{\Omega}$. Then 
\begin{eqnarray*}
\int_{T\left( \triangle \right) }d\tilde{\nu}\left( X\right)
&=&\int_{T\left( \triangle \right) }\tilde{\delta}\left( X\right)
\sup_{Y,Z\in B_{\alpha \tilde{\delta}\left( X\right) }\left( X\right)
}\left| g\left( Y\right) -g\left( Z\right) \right| ^{2}dX \\
&\le &C\int_{T\left( \triangle \right) }\tilde{\delta}\left( X\right) \frac{C%
}{t^{2}}\frac{\sigma \left( \triangle _{X}\right) }{t^{n}}dX \\
&\le &C\frac{\sigma \left( \triangle _{X_{0}}\right) }{t_{0}^{n+1}}\left|
T\left( \triangle \right) \right| \\
&\le &C\frac{d_{0}^{n+1}}{t_{0}^{n+1}}\sigma \left( \triangle
_{X_{0}}\right) \le C\tilde{\mu}\left( \triangle \right) .
\end{eqnarray*}
In case 2, when $d_{0}=\limfunc{diam}T\left( \triangle \right) >\frac{t_{0}}{%
10}$, let $Q_{0}=Q_{cd_{0}}\left( x_{0},0\right) $ be the cube in $\Bbb{R}%
^{n+1}$ centered at $\left( x_{0},0\right) $, with faces parallel to the
coordinate axes and side-length $cd_{0}$. For $c$ big enough depending only
on the Lipschitz character of $\tilde{\Omega}$, we have $T\left( \triangle
\right) \subset \tilde{T}=Q_{0}\bigcap \Bbb{R}_{+}^{n+1}$. Then, since $%
\tilde{\delta}\left( X\right) \le \delta \left( X\right) $ for all $X\in 
\tilde{\Omega}$ and $d\nu \left( X\right) \in \frak{C}\left( d\sigma ,\Omega
\right) $, 
\begin{eqnarray}
\int_{T\left( \triangle \right) }d\tilde{\nu}\left( X\right)
&=&\int_{T\left( \triangle \right) }\tilde{\delta}\left( X\right) \limfunc{%
osc}\nolimits_{\alpha \tilde{\delta}\left( X\right) }g\left( X\right) ^{2}dX
\notag \\
&\le &\int_{\tilde{T}}\delta \left( X\right) \limfunc{osc}\nolimits_{\alpha
\delta \left( X\right) }g\left( X\right) ^{2}dX  \notag \\
&=&\int_{\tilde{T}}d\nu \left( X\right)  \notag \\
&\le &C\mu \left( Q_{0}\bigcap \Bbb{R}^{n}\times \left\{ 0\right\} \right) 
\notag \\
&\le &C\mu \left( Q_{d_{0}}\left( x_{0},0\right) \bigcap \Bbb{R}^{n}\times
\left\{ 0\right\} \right) ,  \label{muti1}
\end{eqnarray}
where we used the doubling property of $\sigma $. Let $\left\{ \mathcal{Q}%
_{i}\right\} _{i=1}^{\infty }$ be a Whitney decomposition of $\Bbb{R}%
_{+}^{n+1}$ into cubes, i.e. for each $i$, $\limfunc{diam}\left( \mathcal{Q}%
_{i}\right) \approx \limfunc{dist}\left( \mathcal{Q}_{i},\partial \Bbb{R}%
_{+}^{n+1}\right) $, for different indexes $\iota $ and $j$, $\mathcal{Q}%
_{i}\bigcap \mathcal{Q}_{j}$ has no interior , and $\Bbb{R}%
_{+}^{n+1}=\bigcup_{i=1}^{\infty }\mathcal{Q}_{i}$. Let $\left\{ \mathcal{P}%
_{j}\right\} _{j=1}^{\infty }$ $\subset \left\{ \mathcal{Q}_{i}\right\}
_{i=1}^{\infty }$ be the collection of cubes such that $\mathcal{P}%
_{j}\bigcap \triangle \neq \emptyset $, and let $X_{j}$ be an arbitrary
point in $\mathcal{P}_{j}\bigcap \triangle $. Then from the definition of $%
\tilde{\mu}$ and since $\limfunc{diam}\left( \mathcal{Q}_{i}\right) \approx 
\limfunc{dist}\left( \mathcal{Q}_{i},\partial \Bbb{R}_{+}^{n+1}\right) $, 
\begin{eqnarray*}
\tilde{\mu}\left( \mathcal{P}_{j}\bigcap \triangle \right) &=&\int_{\mathcal{%
P}_{j}\bigcap \triangle }M\left( X\right) d\tilde{\sigma}\left( X\right)
\approx \frac{\sigma \left( \triangle _{X_{j}}\right) }{\sigma \left(
\triangle _{X_{j}}\right) }\int_{\mathcal{P}_{j}\bigcap \triangle }d\tilde{%
\sigma}\left( X\right) \\
&\approx &C\frac{\sigma \left( \triangle _{X_{j}}\right) }{\limfunc{diam}%
\left( \mathcal{P}_{i}\right) ^{n}}\int_{\mathcal{P}_{j}\bigcap \triangle }d%
\tilde{\sigma}\left( X\right) \approx C\mu \left( \triangle _{X_{j}}\right) .
\end{eqnarray*}
Then 
\begin{eqnarray}
\tilde{\mu}\left( \triangle \right) &=&\sigma \left( \triangle \bigcap
\Omega \right) +\sum_{j=1}^{\infty }\tilde{\mu}\left( \mathcal{P}_{j}\bigcap
\triangle \right) \approx \sigma \left( \triangle \bigcap \Omega \right)
+\sum_{j=1}^{\infty }\sigma \left( \triangle _{X_{j}}\right)  \label{muti3}
\\
&\ge &C\mu \left( Q_{d_{0}}\left( x_{0},0\right) \bigcap \Bbb{R}^{n}\times
\left\{ 0\right\} \right) ,  \notag  \label{muti2}
\end{eqnarray}
the last inequality follows from a simple geometrical argument. From this
and (\ref{muti1}) we have $\int_{T\left( \triangle \right) }d\tilde{\nu}%
\left( X\right) \le C\tilde{\mu}\left( \triangle \right) $ as wanted. Also,
from (\ref{muti3}) we also have $\tilde{\mu}\left( \triangle \right) \approx 
\tilde{\mu}\left( 2\triangle \right) $ in this case. 
\endproof%

The following lemma estates the local character of the regularity of the
harmonic measure.

\begin{lemma}
\label{main}Let $\mathcal{L}$ be an elliptic operator in divergence form or
nondivergence form with drift $\mathbf{b}$ in a Lipschitz domain $\Omega $;
i.e. $\mathcal{L}=\mathbf{A}\cdot \nabla +\mathbf{b}\cdot \nabla $ or $%
\mathcal{L=}\mathbf{A}\cdot \nabla ^{2}+\mathbf{b}\cdot \nabla $ where $%
\mathbf{A}$ satisfies the ellipticity condition (\ref{ellipticity}). Suppose
that $\mathbf{b\,}$is locally bounded in $\Omega $ and that for a doubling
measure $\sigma $on $\partial \Omega $ it satisfies 
\begin{equation}
\delta \left( X\right) \limfunc{osc}\nolimits_{\frac{\delta \left( X\right) 
}{2\sqrt{n}}}\mathbf{b}\left( X\right) ^{2}\,\in \frak{C}\left( d\sigma
,\Omega \right) \text{.}  \label{carlb}
\end{equation}
Then $d\omega _{\mathcal{L}}\in A_{\infty }$ \emph{if and only if} there
exists a finite collection of Lipschitz domains $\left\{ \Omega \right\}
_{i=1}^{N}$ and compact sets $K_{i}\subset \subset \partial \Omega
_{i}\bigcap \partial \Omega $ such that $\bigcup_{i=1}^{n}\Omega _{i}\subset
\Omega $, $\partial \Omega \subset \bigcup_{i=1}^{n}K_{i}$, and 
\begin{equation*}
d\omega _{\mathcal{L}_{i}|K_{i}}\in A_{\infty }\left( d\sigma \right)
,\qquad i=1,\dots ,N.
\end{equation*}
where $\mathcal{L}_{i}$ denotes the restriction of $\mathcal{L}$ to the
subdomain $\Omega _{i}$.
\end{lemma}

\proof%
If $\mathbf{b}=0$, the result is an immediate consequence of the ``main
lemma'' in \cite{DJK} for the divergence case and the analog to the main
lemma in the nondivergence case, contained in \cite{EK}. The case $\mathbf{b}%
\neq 0$ then follows from Theorem \ref{HL217} for the divergence case and
Theorem \ref{CRth2b} from next section for the nondivergence case. Indeed,
by the mentioned theorems, if $\mathcal{L}$ is the operator with drift $%
\mathbf{b}$ and $\mathcal{L}_{0}$ is the operator with the same second order
coefficients but \emph{without} a drift term, then $d\omega _{\mathcal{L}%
}\in A_{\infty }\left( d\sigma \right) \Leftrightarrow d\omega _{\mathcal{L}%
_{0}}\in A_{\infty }\left( d\sigma \right) $. On the other hand, by Lemma 
\ref{carlesin} with $g\left( X\right) =\mathbf{b}\left( X\right) $, the
restriction of $\mathbf{b}$ to any Lipschitz subdomain $\Omega ^{\prime
}\subset \Omega $ also satisfies (\ref{carlb}) in $\Omega ^{\prime }$.
Hence, by Theorems \ref{HL217} and \ref{CRth2b}, for any doubling measure $%
d\sigma ^{\prime }$ on $\partial \Omega ^{\prime }$, we have $d\omega _{%
\mathcal{L}|\Omega ^{\prime }}\in A_{\infty }\left( d\sigma ^{\prime
}\right) \Leftrightarrow d\omega _{\mathcal{L}_{0}|\Omega ^{\prime }}\in
A_{\infty }\left( d\sigma ^{\prime }\right) $. Now, for $\Omega _{i}$, $%
K_{i} $ as in the statement of Lemma \ref{main}, $d\sigma |_{K_{i}}$, the
restriction of $d\sigma $ to the compact set $K_{i}$, can be extended to a
doubling measure $d\sigma _{i}$ on $\partial \Omega _{i}$ with the same
doubling constant. Then $d\omega _{\mathcal{L}|\Omega _{i}}\in A_{\infty
}\left( d\sigma ^{i}\right) \Leftrightarrow d\omega _{\mathcal{L}_{0}|\Omega
_{i}}\in A_{\infty }\left( d\sigma ^{i}\right) $, which implies that 
\begin{equation*}
d\omega _{\mathcal{L}|\Omega _{i}}|_{K_{i}}\in A_{\infty }\left( d\sigma
|_{K_{i}}\right) \Leftrightarrow d\omega _{\mathcal{L}_{0}|\Omega
_{i}}|_{K_{i}}\in A_{\infty }\left( d\sigma |_{K_{i}}\right) ,
\end{equation*}
where $d\omega _{\mathcal{L}|\Omega _{i}}|_{K_{i}}$ denotes the restriction
to $K_{i}$ of the harmonic measure of $\mathcal{L}$ in $\Omega _{i}$, with a
similar definition for $d\omega _{\mathcal{L}_{0}|\Omega _{i}}|_{K_{i}}$.
This shows that Lemma \ref{main} in the case $\mathbf{b}\neq 0$ follows from
the case $\mathbf{b}=0$. 
\endproof%

\section{Proofs of the Theorems\label{proofs}}

The proof of Theorem \ref{CRth2} relies on the following special case.

\begin{theorem}
\label{CRth2b} Let $\mathcal{\tilde{L}}_{N,0}=\mathbf{A}\cdot \nabla ^{2}$
and $\mathcal{\tilde{L}}_{N,1}=\mathbf{A}\nabla ^{2}+\mathbf{b}\cdot \nabla $
where $\mathbf{A}=\left( A_{ij}\right) _{i,j=1}^{n+1}$ and $\mathbf{b}%
=\left( b_{j}\right) _{j=1}^{n+1}$ are bounded, measurable coefficients and $%
\mathbf{A}$ satisfy the ellipticity condition (\ref{ellipticity}) in a
Lipschitz domain $\Omega $. Suppose that $\mathcal{CD}$ holds for $\mathcal{%
\tilde{L}}_{N,\ell }$, $\ell =0,1$ in $\Omega $ and that 
\begin{equation}
\delta \left( X\right) \sup_{Y\in Q_{\frac{\delta \left( X\right) }{2\sqrt{n}%
}}\left( X\right) }\left| \mathbf{b}\left( Y\right) \right| ^{2}dX\in \frak{C%
}\left( d\sigma ,\Omega \right) .  \label{CarlCRb}
\end{equation}
Then $d\omega _{\mathcal{\tilde{L}}_{N,0}}\in A_{\infty }\Rightarrow d\omega
_{\mathcal{\tilde{L}}_{N,1}}\in A_{\infty }$.
\end{theorem}

We defer the proof of this result (which contains the main substance of
Theorem \ref{CRth2}) to next section. Now we obtain Theorem \ref{CRth2} from
Theorem \ref{CRth2b}.

\proof[Proof of Theorem~\ref{CRth2}]%
Let $\mathcal{L}_{N,0}=\mathbf{A}_{0}\cdot \nabla ^{2}+\mathbf{b}_{0}\cdot
\nabla $ and $\mathcal{L}_{N,1}=\mathbf{A}_{1}\nabla ^{2}+\mathbf{b}%
_{1}\cdot \nabla $ where $\mathbf{A}_{\ell },$ $\mathbf{b}_{\ell }$ and $%
\mathcal{L}_{N,\ell }$ satisfy the hypotheses of Theorem \ref{CRth2} for $%
\ell =0,1$. Then if $d\omega _{\mathcal{L}_{N,0}}\in A_{\infty }$, by
Theorem \ref{CRth2b} it follows that $d\omega _{\mathcal{L}_{N,0}^{*}}\in
A_{\infty }$ where $\mathcal{L}_{N,0}^{*}=\mathbf{A}_{0}\cdot \nabla ^{2}$.
Indeed, since obviously $d\omega _{\mathcal{L}_{N,0}^{*}}\in A_{\infty
}\left( d\omega _{\mathcal{L}_{N,0}^{*}}\right) $, by Theorem \ref{CRth2b}
we have that $d\omega _{\mathcal{L}_{N,0}}\in A_{\infty }\left( d\omega _{%
\mathcal{L}_{N,0}^{*}}\right) \bigcap A_{\infty }$; and hence $d\omega _{%
\mathcal{L}_{N,0}^{*}}\in A_{\infty }$. By the result in \cite{CR} (Theorem 
\ref{FKP0} in the Introduction) and Theorem \ref{CRth2b}, we have that if $%
\mathcal{L}_{N,1}^{*}=\mathbf{A}_{1}\cdot \nabla ^{2}$, then 
\begin{equation*}
d\omega _{\mathcal{L}_{N,0}^{*}}\in A_{\infty }\Rightarrow d\omega _{%
\mathcal{L}_{N,1}^{*}}A_{\infty }\Rightarrow d\omega _{\mathcal{L}_{N,1}}\in
A_{\infty },
\end{equation*}
as wanted.%
\endproof%

Theorem \ref{CRth1} will follow once we are able to apply similar
techniques. This will be accomplished with the aid of the approximation
provided by Lemma \ref{grcarl} in the special case $\Omega =\Bbb{R}%
_{+}^{n+1} $, and the localization given by Lemma \ref{main}.

\proof[Proof of Theorem~\ref{CRth1}]%
Given $P_{0}\in \partial \Omega $ let $X=\left( x,t\right) \,$be a
coordinate system such that $P_{0}=\left( x_{0},t_{0}\right) \in \partial
\Omega $ and there exists a Lipschitz function $\psi :\Bbb{R}^{n}\rightarrow 
\Bbb{R}$ defining a local coordinate system of $\Omega $ in a neighborhood
of $P_{0}$. That is, for some $r_{0}=r_{0}\left( \Omega \right) >0$ we have 
\begin{eqnarray}
\partial \Omega \bigcap \left\{ \left| x-x_{0}\right| <r_{0}\right\} \times 
\Bbb{R} &=&\left\{ \left( x,\psi \left( x\right) \right) :\left|
x-x_{0}\right| <r_{0}\right\}  \notag \\
\Omega ^{\prime } &=&\left\{ \left( x,t\right) :\left| x-x_{0}\right|
<r_{0},\,\psi \left( x\right) <t<\psi \left( x\right) +r_{0}\right\} \subset
\Omega .  \label{op}
\end{eqnarray}
Let $\eta _{s}\left( y\right) =s^{-n}\eta \left( \frac{y}{s}\right) $, where 
$\eta $ is an even $\mathcal{C}^{\infty }$ approximate identity in $\Bbb{R}%
^{n}$ supported in $\left\{ \left| y\right| \le \frac{1}{2}\right\} $. Set $%
\rho \left( y,s\right) =\left( y,c_{0}s+F\left( y,s\right) \right) $ with $%
F\left( y,s\right) =\eta _{s}*\psi \left( y\right) =\int_{\Bbb{R}^{n}}\eta
_{s}\left( y-z\right) \psi \left( z\right) dz$. We have 
\begin{equation*}
\nabla _{Y}\rho =\left( 
\begin{array}{ll}
I & \nabla _{y}F \\ 
0 & c_{0}+\frac{\partial F}{\partial s}
\end{array}
\right) ,
\end{equation*}
Since $\left\| \frac{\partial F}{\partial s}\right\| _{\infty }\le
C_{n}\left\| \nabla _{y}\psi \right\| _{\infty }$, where $C_{n}=n\int \eta
\left( y\right) \left| y\right| d\sigma $, (note that, for appropriate $\eta 
$, $C_{n}$ is a universal constant), taking $c_{0}=1+C_{n}\left\| \nabla
_{y}\psi \right\| _{\infty }$, $\rho $ is a 1-1 map of $\Bbb{R}_{+}^{n+1}$
onto $\left\{ \left( x,t\right) :t>\psi \left( x\right) \right\} $,
moreover, $\rho $ is bi-Lipschitz and $1\le \left| \det \nabla _{Y}\rho
\right| \le 1+2C_{n}\left\| \nabla _{y}\psi \right\| _{\infty }$. This
transformation gives rise to the Dahlberg-Kenig-Stein adapted distance
function. For $\alpha >0$, let $\Phi _{\alpha }=\left\{ y:\left|
y-x_{0}\right| <\alpha \right\} \times \left( 0,\frac{\alpha }{c_{0}}\right) 
$, and $\Omega _{\alpha }=\rho \left( \Phi _{\alpha }\right) $. For $\alpha
=\alpha _{0}$ small enough depending only on $\left\| \nabla _{y}\psi
\right\| _{\infty }$ and $n$, we have $\Omega _{\alpha _{0}}\subset \Omega
^{\prime }$, where $\Omega ^{\prime }$ is given by\newline
(\ref{op}). Moreover, $\Omega _{\alpha _{0}}\,$is a Lipschitz domain with
Lipschitz constant depending only on the constant of $\Omega $.

We will first consider the divergence case, suppose that $\mathcal{L}_{D}=%
\limfunc{div}\mathbf{A}\nabla +\mathbf{b}\cdot \nabla $ is a uniformly
elliptic operators in divergence form with bounded measurable coefficient
matrix $\mathbf{A}$ and drift vector $\mathbf{b}$ satisfying (\ref{CarlCRa}%
), i.e. 
\begin{equation}
\left\{ \delta ^{-1}\left( X\right) \limfunc{osc}\nolimits_{Q_{\frac{\delta
\left( X\right) }{2\sqrt{n}}}\left( X\right) }\mathbf{A}\left( X\right)
^{2}+\delta \left( X\right) \limfunc{osc}\nolimits_{Q_{\frac{\delta \left(
X\right) }{2\sqrt{n}}}\left( X\right) }\mathbf{b}\left( X\right)
^{2}\right\} dX\in \frak{C}\left( d\sigma ,\Omega \right) .  \label{CarlCRa1}
\end{equation}

Since $P_{0}$ is an arbitrary point on $\partial \Omega $ and $\partial
\Omega $ is compact, by Lemma \ref{main}, to prove Theorem \ref{CRth1} it is
enough to prove that if $\omega $ is the harmonic measure for $\mathcal{L}%
_{D}\,$in $\Omega _{\alpha _{0}}$, then 
\begin{equation}
\omega |_{K}\in A_{\infty }\left( d\sigma |_{K}\right) ,  \label{enough}
\end{equation}
where $K\subset \partial \Omega _{\alpha _{0}}$ is the compact set given by 
\begin{equation}
K=\rho \left\{ Y=\left( y,s\right) :\left| y-x_{0}\right| \le \frac{\alpha
_{0}}{6},s=0\right\} \subset \partial \Omega _{\alpha _{0}}.  \label{K}
\end{equation}
For simplicity, we will write $\Omega =\Omega _{\alpha _{0}}$and $\Phi =\Phi
_{\alpha _{0}}$. Since $\psi $ is Lipschitz, it follows that the
transformation $\rho :\Phi \rightarrow \Omega $ can be extended to a
homeomorphism from $\overline{\Phi }$ to $\overline{\Omega }$ and such that
the restriction of $\rho $ to $\partial \Phi $ is a bi-Lipschitz
homeomorphism from $\partial \Omega $ to $\partial \Omega $. Indeed, since $%
\left\{ \eta _{s}\right\} _{s>0}$ is a smooth approximation of the identity,
it easily follows that $\rho $ restricted to $\partial \Phi $ is given by 
\begin{equation}
\rho \left( Y\right) =\rho \left( y,s\right) =\left\{ 
\begin{array}{lll}
\left( y,\psi \left( y\right) \right) , &  & s=0 \\ 
&  &  \\ 
\left( y,c_{0}s+F\left( y,s\right) \right) , &  & 0<s<\frac{\alpha _{0}}{%
c_{0}},\quad \left| y-x_{0}\right| =\alpha _{0} \\ 
&  &  \\ 
\left( y,c_{0}\frac{\alpha _{0}}{c_{0}}+F\left( y,\frac{\alpha _{0}}{c_{0}}%
\right) \right) &  & s=\frac{\alpha _{0}}{c_{0}}
\end{array}
\right.  \label{rY}
\end{equation}
$\qquad $whenever $Y\in \partial \Phi $. If $\tilde{\delta}\left( Y\right) =%
\limfunc{dist}\left( Y,\partial \Phi \right) $, then for some constant $%
\tilde{C}$ depending only on $\left\| \nabla _{y}\psi \right\| _{\infty }$
and $n$, the following estimate holds for the distance functions $\delta $
and $\tilde{\delta}$: 
\begin{equation}
\tilde{C}^{-1}\tilde{\delta}\left( Y\right) \le \delta \left( \rho \left(
Y\right) \right) \le \tilde{C}\delta \left( Y\right) ,\qquad Y\in \Phi \text{%
.}  \label{deltas}
\end{equation}
To see this, let $Y_{0}\in \Phi $ and let $\tilde{\delta}_{0}=\delta \left(
Y_{0}\right) $. Now, $\delta \left( \rho \left( Y_{0}\right) \right) =%
\limfunc{dist}\left( \rho \left( Y_{0}\right) ,\partial \Omega \right)
=\left| \rho \left( Y_{0}\right) -X^{\prime }\right| $ for some $X^{\prime
}\in \partial \Omega $. Let $Y^{\prime }=\rho ^{-1}\left( X^{\prime }\right) 
$, and $X_{0}=\rho \left( Y_{0}\right) $, thus, $Y^{\prime }\in \partial
\Phi $ and since $\rho ^{-1}$ is Lipschitz in $\overline{\Phi }$, we have 
\begin{equation*}
\tilde{\delta}\left( Y_{0}\right) \le \left| Y_{0}-Y^{\prime }\right|
=\left| \rho ^{-1}\left( X_{0}\right) -\rho ^{-1}\left( X^{\prime }\right)
\right| \le C\left| X_{0}-X^{\prime }\right| =C\delta \left( \rho \left(
Y_{0}\right) \right) ,
\end{equation*}
with $C=C\left( n,\left\| \nabla _{y}\psi \right\| _{\infty }\right) $. The
other inequality in (\ref{deltas}) follows in a similar manner.

The Lebesgue measure $\sigma $ on $\partial \Omega $ induces a doubling
measure $\tilde{\mu}$ on $\partial \Phi $ by the relation 
\begin{equation*}
\tilde{\mu}\left( E\right) =\sigma \left( \rho \left( E\right) \right) \text{%
,\quad for any Borel set }E\subset \partial \Phi .
\end{equation*}
Let $\tilde{\sigma}$ be the Lebesgue measure on $\partial \Phi $, then from
the definition of $\tilde{\mu}$ and the fact that $\rho $ is bi-Lipschitz it
easily follows that $\frac{d\tilde{\mu}}{d\tilde{\sigma}}\approx 1$. Hence
we can replace $\tilde{\mu}$ by $\tilde{\sigma}$ in our calculations. Now,
if $u\left( x,t\right) $ is a solution of $\mathcal{L}_{D}u=\limfunc{div}_{X}%
\mathbf{A}\nabla _{X}u+\mathbf{b}\cdot \nabla _{X}u=0$ in $\Omega \,$, then $%
v\left( y,s\right) =u\left( \rho \left( y,s\right) \right) $, defined in $%
\Phi ,$ is a solution of $\mathcal{\tilde{L}}_{D}v=\limfunc{div}_{Y}\mathbf{%
\tilde{A}}\nabla _{Y}v+\mathbf{\tilde{b}}\cdot \nabla _{Y}v=0$, where 
\begin{eqnarray}
\mathbf{\tilde{A}}\left( Y\right) &=&\left( \left( \nabla _{Y}\rho \right)
^{-1}\left( Y\right) \right) ^{t}\mathbf{A}\left( \rho \left( Y\right)
\right) \left( \nabla _{Y}\rho \right) ^{-1}\left( Y\right) \det \left(
\nabla _{Y}\rho \right) \left( Y\right) ,  \label{Atilde} \\
\mathbf{\tilde{b}}\left( Y\right) &=&\mathbf{b}\left( \rho \left( Y\right)
\right) \mathbf{\cdot }\,\left( \nabla _{Y}\rho \right) ^{-1}\left( Y\right)
\det \left( \nabla _{Y}\rho \right) \left( Y\right) .  \label{btilde}
\end{eqnarray}
We claim that $\mathbf{\tilde{A}}$ and $\mathbf{\tilde{b}}$ satisfy 
\begin{equation}
\left\{ \tilde{\delta}\left( Y\right) ^{-1}\limfunc{osc}\nolimits_{Q_{\frac{%
\delta \left( Y\right) }{2\sqrt{n}}}\left( Y\right) }\mathbf{\tilde{A}}%
\left( Y\right) ^{2}+\tilde{\delta}\left( Y\right) \limfunc{osc}%
\nolimits_{Q_{\frac{\delta \left( Y\right) }{2\sqrt{n}}}\left( Y\right) }%
\mathbf{b}\left( Y\right) ^{2}\right\} dY\in \frak{C}\left( d\tilde{\sigma}%
,\Phi \right) ,  \label{CarlCRa1T}
\end{equation}
Let $Y_{0}\in \Phi $ let $Y_{1}$, $Y_{2}$ such that $\left|
Y_{1}-Y_{0}\right| \le \frac{1}{2}\delta _{0}$ and $\left|
Y_{2}-Y_{0}\right| \le \frac{1}{2}\delta _{0}$, where $\delta _{0}=\tilde{%
\delta}\left( Y_{0}\right) $, then by (\ref{Atilde}) 
\begin{eqnarray}
&&\frac{\left| \mathbf{\tilde{A}}\left( Y_{1}\right) -\mathbf{\tilde{A}}%
\left( Y_{2}\right) \right| ^{2}}{\tilde{\delta}_{0}}  \notag \\
&=&\tilde{\delta}_{0}^{-1}\left| \left( \left( \nabla _{Y}\rho \right)
^{-1}\left( Y_{1}\right) \right) ^{t}\mathbf{A}\left( \rho \left(
Y_{1}\right) \right) \left( \nabla _{Y}\rho \right) ^{-1}\left( Y_{1}\right)
\det \left( \nabla _{Y}\rho \right) \left( Y_{1}\right) \right.  \notag \\
&&\left. -\left( \left( \nabla _{Y}\rho \right) ^{-1}\left( Y_{2}\right)
\right) ^{t}\mathbf{A}\left( \rho \left( Y_{2}\right) \right) \left( \nabla
_{Y}\rho \right) ^{-1}\left( Y_{2}\right) \det \left( \nabla _{Y}\rho
\right) \left( Y_{2}\right) \right| ^{2}  \notag \\
&\le &C\frac{\left| \left( \left( \nabla _{Y}\rho \right) ^{-1}\left(
Y_{1}\right) \right) ^{t}-\left( \left( \nabla _{Y}\rho \right) ^{-1}\left(
Y_{2}\right) \right) ^{t}\right| ^{2}}{\tilde{\delta}_{0}}  \label{dleltaA}
\\
&&\cdot \left| \mathbf{A}\left( \rho \left( Y_{1}\right) \right) \left(
\nabla _{Y}\rho \right) ^{-1}\left( Y_{1}\right) \det \left( \nabla _{Y}\rho
\right) \left( Y_{1}\right) \right| ^{2}  \notag \\
&&+C\frac{\left| \left( \nabla _{Y}\rho \right) ^{-1}\left( Y_{1}\right)
-\left( \nabla _{Y}\rho \right) ^{-1}\left( Y_{2}\right) \right| ^{2}}{%
\tilde{\delta}_{0}}  \notag \\
&&\cdot \left| \left( \left( \nabla _{Y}\rho \right) ^{-1}\left(
Y_{2}\right) \right) ^{t}\mathbf{A}\left( \rho \left( Y_{1}\right) \right)
\right| ^{2}\left| \det \left( \nabla _{Y}\rho \right) \left( Y_{1}\right)
\right| ^{2}  \notag \\
&&+C\frac{\left| \det \left( \nabla _{Y}\rho \right) \left( Y_{1}\right)
-\det \left( \nabla _{Y}\rho \right) \left( Y_{2}\right) \right| ^{2}}{%
\tilde{\delta}_{0}}  \notag \\
&&\cdot \left| \left( \left( \nabla _{Y}\rho \right) ^{-1}\left(
Y_{2}\right) \right) ^{t}\mathbf{A}\left( \rho \left( Y_{1}\right) \right)
\left( \nabla _{Y}\rho \right) ^{-1}\left( Y_{2}\right) \right| ^{2}  \notag
\\
&&+C\frac{\left| \mathbf{A}\left( \rho \left( Y_{1}\right) \right) -\mathbf{A%
}\left( \rho \left( Y_{2}\right) \right) \right| ^{2}}{\tilde{\delta}_{0}} 
\notag \\
&&\cdot \left| \left( \left( \nabla _{Y}\rho \right) ^{-1}\left(
Y_{2}\right) \right) ^{t}\right| ^{2}\left| \left( \nabla _{Y}\rho \right)
^{-1}\left( Y_{2}\right) \det \left( \nabla _{Y}\rho \right) \left(
Y_{2}\right) \right| ^{2}.  \notag
\end{eqnarray}
and similarly, by (\ref{btilde}) 
\begin{eqnarray}
&&\tilde{\delta}_{0}\left| \mathbf{b}\left( Y_{1}\right) -\mathbf{b}\left(
Y_{2}\right) \right| ^{2}  \notag \\
&=&\tilde{\delta}_{0}\left| \mathbf{b}\left( \rho \left( Y_{1}\right)
\right) \mathbf{\cdot }\,\left( \nabla _{Y}\rho \right) ^{-1}\left(
Y_{1}\right) \det \left( \nabla _{Y}\rho \right) \left( Y_{1}\right) \right.
\notag \\
&&\left. -\mathbf{b}\left( \rho \left( Y_{2}\right) \right) \mathbf{\cdot }%
\,\left( \nabla _{Y}\rho \right) ^{-1}\left( Y_{2}\right) \det \left( \nabla
_{Y}\rho \right) \left( Y_{2}\right) \right| ^{2}  \notag \\
&\le &C\tilde{\delta}_{0}\left| \mathbf{b}\left( \rho \left( Y_{1}\right)
\right) -\mathbf{b}\left( \rho \left( Y_{2}\right) \right) \right|
^{2}\left| \mathbf{\cdot }\left( \nabla _{Y}\rho \right) ^{-1}\left(
Y_{1}\right) \det \left( \nabla _{Y}\rho \right) \left( Y_{1}\right) \right|
^{2}  \label{dlelta} \\
&&+C\tilde{\delta}_{0}\left| \mathbf{b}\left( \rho \left( Y_{2}\right)
\right) \mathbf{\cdot }\,\left( \nabla _{Y}\rho \right) ^{-1}\left(
Y_{1}\right) \right| ^{2}\left| \det \left( \nabla _{Y}\rho \right) \left(
Y_{1}\right) -\det \left( \nabla _{Y}\rho \right) \left( Y_{1}\right)
\right| ^{2}  \notag \\
&&+C\tilde{\delta}_{0}\left| \mathbf{b}\left( \rho \left( Y_{2}\right)
\right) \right| ^{2}\,\left| \left( \nabla _{Y}\rho \right) ^{-1}\left(
Y_{1}\right) -\left( \nabla _{Y}\rho \right) ^{-1}\left( Y_{2}\right)
\right| ^{2}\left| \det \left( \nabla _{Y}\rho \right) \left( Y_{2}\right)
\right| ^{2}.  \notag
\end{eqnarray}
We will use the following fact:

\begin{lemma}
\label{rho}For $\Omega $, $\Phi $, $\sigma $, $\tilde{\mu}$, $\delta $, $%
\tilde{\delta}$ and $\rho $ as above, the functions 
\begin{equation*}
r_{1}\left( Y\right) =\frac{1}{\tilde{\delta}\left( Y\right) }\left| 
\limfunc{osc}\nolimits_{Q_{\frac{\tilde{\delta}\left( Y\right) }{2\sqrt{n}}%
}\left( Y\right) }\left( \nabla _{Y}\rho \right) ^{-1}\left( Y\right)
\right| ^{2}
\end{equation*}
and 
\begin{equation*}
r_{2}\left( Y\right) =\frac{1}{\tilde{\delta}\left( Y\right) }\left| 
\limfunc{osc}\nolimits_{Q_{\frac{\tilde{\delta}\left( Y\right) }{2\sqrt{n}}%
}\left( Y\right) }\left( \det \left( \nabla _{Y}\rho \right) \left( Y\right)
\right) \right| ^{2}
\end{equation*}
defined in $\Phi $, satisfy $\left[ r_{1}\left( Y\right) +r_{2}\left(
Y\right) \right] dY$ $\in \frak{C}\left( d\tilde{\sigma},\Phi \right) $,
where $\tilde{\sigma}\,$is the Lebesgue measure on $\partial \Phi $.
\end{lemma}

The lemma follows from the fact that $\tilde{\delta}\left( Y\right) \left|
\nabla ^{2}\rho \right| ^{2}dY\in \frak{C}\left( d\tilde{\sigma},\Phi
\right) $. This property is discussed in \cite{KP}, and it can be obtained
as an application of the characterization of $A_{\infty }$ in terms of
Carleson measures given in \cite{FKP}.

Then, (\ref{CarlCRa1T}) follows by applying (\ref{CarlCRa1}), (\ref{deltas})
and Lemma \ref{rho} to the expressions (\ref{dleltaA}) and (\ref{dlelta}).

We recall that $\Phi =\Phi _{\alpha _{0}}$ where $\Phi _{\alpha }$ is given
by $\Phi _{\alpha }=\left\{ y:\left| y-x_{0}\right| <\alpha \right\} \times
\left( 0,\frac{\alpha }{c_{0}}\right) $. Denote by $\Phi _{\alpha }^{\pm
}=\left\{ y:\left| y-x_{0}\right| <\alpha \right\} \times \left( -\frac{%
\alpha }{c_{0}},\frac{\alpha }{c_{0}}\right) $ and let $\nu \left( Y\right)
\in \mathcal{C}_{0}^{\infty }\left( \Phi _{\alpha _{0}}^{\pm }\right) $ such
that $0\le \nu \le 1$, $\nu \equiv 1$ in $\Phi _{\frac{\alpha _{0}}{3}}^{\pm
}$ and $\nu \equiv 0$ in $\Phi _{\alpha _{0}}^{\pm }\backslash \Phi _{\frac{%
2\alpha _{0}}{3}}^{\pm }$. For $Y\in \Bbb{R}_{+}^{n+1}$, let 
\begin{equation*}
\mathbf{\tilde{A}}^{*}\left( Y\right) =\nu \left( Y\right) \mathbf{\tilde{A}}%
\left( Y\right) +\left( 1-\nu \left( Y\right) \right) I,
\end{equation*}
where $I$ is the $\left( n+1\right) \times \left( n+1\right) $ identity
matrix. It follows that $\mathbf{\tilde{A}}^{*}\left( Y\right) $ is an
elliptic matrix function, with the same ellipticity constants as $\mathbf{%
\tilde{A}}$. The measure $\tilde{\sigma}$ extends trivially from $\partial
\Phi \bigcap \partial \Bbb{R}_{+}^{n+1}$ to $\partial \Bbb{R}_{+}^{n+1}$, we
dub this extension (which is just the Euclidean measure) $d\tilde{\sigma}%
^{*} $. With this definitions, because of (\ref{CarlCRa1T}), for $Y=\left(
y,t\right) \in \Bbb{R}_{+}^{n+1}$ , $\mathbf{\tilde{A}}^{*}$ satisfies 
\begin{equation*}
\frac{\limfunc{osc}\nolimits_{Q_{\frac{t}{2\sqrt{n}}}\left( Y\right) }%
\mathbf{\tilde{A}}^{*}\left( Y\right) ^{2}}{t}\in \frak{C}\left( d\tilde{%
\sigma}^{*},\Bbb{R}_{+}^{n+1}\right) ,
\end{equation*}
(i.e.: is a Carleson measure in $\Bbb{R}_{+}^{n+1}$ with respect to $\tilde{%
\sigma}^{*}$). Where $Q_{\gamma }\left( Y\right) $ is the cube centered at $%
Y $ with faces parallel to the coordinate axes and sidelength $\gamma t$.
Therefore, $\mathbf{\tilde{A}}^{*}\left( Y\right) $ satisfies the hypotheses
of Lemma \ref{grcarl}, and the matrix function $\mathbf{\tilde{A}}%
^{**}\left( Y\right) $ defined by 
\begin{equation}
\mathbf{\tilde{A}}^{**}\left( Y\right) =\frac{1}{\left| Q_{0}\left( Y\right)
\right| }\int_{Q_{0}\left( Y\right) }\mathbf{\hat{A}}\left( Z\right) dZ,
\label{Ass}
\end{equation}
where $Q_{0}\left( Y\right) =Q_{\frac{\delta \left( Y\right) }{6\sqrt{n}}%
}\left( Y\right) $ and $\mathbf{\hat{A}}\left( Z\right) =\tfrac{1}{\left|
Q_{0}\left( Z\right) \right| }\int_{Q_{0}\left( Z\right) }\mathbf{\tilde{A}}%
^{*}\left( W\right) dW$, satisfies 
\begin{equation}
t\mathcal{E}^{*}\left( Y\right) ^{2}\in \frak{C}\left( d\tilde{\sigma}^{*},%
\Bbb{R}_{+}^{n+1}\right) \quad \text{and}\quad \frac{\mathcal{\tilde{E}}%
^{**}\left( X\right) ^{2}}{t}dX\in \frak{C}\left( d\tilde{\sigma}^{*},\Bbb{R}%
_{+}^{n+1}\right)  \label{Ezz}
\end{equation}
where 
\begin{equation}
\mathcal{\tilde{E}}^{*}\left( Y\right) =\sup_{Z\in Q_{0}\left( Y\right)
}\left| \nabla \mathbf{\tilde{A}}^{**}\left( Z\right) \right| \quad \text{and%
}\quad \mathcal{\tilde{E}}^{**}\left( Y\right) =\sup_{Z\in Q_{0}\left(
Y\right) }\left| \mathbf{\tilde{A}}^{**}\left( Z\right) -\mathbf{\tilde{A}}%
^{*}\left( Z\right) \right| .  \label{Ess}
\end{equation}
Moreover, from the definitions it is easy to check that $\mathbf{\tilde{A}}%
^{**}$ is elliptic with the same ellipticity constants as $\mathbf{\tilde{A}}
$.

Let now $\mathcal{\tilde{L}}_{D}^{*}=\limfunc{div}_{Y}\mathbf{\tilde{A}}%
^{*}\nabla _{Y}+\mathbf{\tilde{b}}\cdot \nabla _{Y}$ and $\mathcal{\tilde{L}}%
_{D}^{**}=\limfunc{div}_{Y}\mathbf{\tilde{A}}^{**}\nabla _{Y}+\mathbf{\tilde{%
b}}\cdot \nabla _{Y}$. By the Carleson measure property of $\mathcal{\tilde{E%
}}^{**}$, and Lemma \ref{carlesin}, $\mathcal{\tilde{L}}_{D}^{*}$ and $%
\mathcal{\tilde{L}}_{D}^{**}$ satisfy the hypotheses of Theorem \ref{HL217}
in $\Phi $ with respect to the measure $\tilde{\mu}$, therefore, if $\tilde{%
\omega}^{*}$ and $\tilde{\omega}^{**}$ denote the harmonic measures of $%
\mathcal{\tilde{L}}_{D}^{*}$ and $\mathcal{\tilde{L}}_{D}^{**}$ in $\Phi $,
respectively, we have that $\tilde{\omega}^{*}\in A_{\infty }\Leftrightarrow 
\tilde{\omega}^{**}\in A_{\infty }$. On the other hand, because of the
Carleson measure property of $\mathcal{\tilde{E}}^{*}$, and Lemma \ref
{carlesin}, $\mathcal{\tilde{L}}_{D}^{**}$ satisfies the hypotheses of
Theorem \ref{KP218} in $\Phi $; and therefore $\tilde{\omega}^{**}\in
A_{\infty }$. From what we just proved it follows that $\tilde{\omega}%
^{*}\in A_{\infty }$.

Let $\mathcal{L}_{D}^{*}$ denote the pull-back of $\mathcal{\tilde{L}}%
_{D}^{*}$ from $\Phi $ to $\Omega $ through the mapping $\rho $. Since the
mapping $\rho :\partial \Phi \rightarrow \partial \Omega \,$is bi-Lipschitz,
and if $\omega ^{*}$ denotes the harmonic measure of $\mathcal{L}_{D}^{*}%
\mathcal{\ }$in $\Omega $, then $\omega ^{*}\in A_{\infty }$. This follows
directly from the definitions of the $A_{\infty }$ class, the harmonic
measure $\omega ^{*}$ and $\mathcal{L}_{D}^{*}$. On the other hand, the
operator $\mathcal{L}_{D}^{*}$ coincides with $\mathcal{L}_{D}=\limfunc{div}%
\nolimits_{X}\mathbf{A}\nabla _{X}+\mathbf{b}\cdot \nabla _{X}$ in $\rho
\left( \Phi _{\frac{\alpha _{0}}{3}}\right) $, hence an application of
Theorem \ref{HL217} and the``main lemma'' in \cite{DJK} (see Lemma \ref{main}%
), implies that $\omega |_{K}\in A_{\infty }\left( d\omega ^{*}|_{K}\right) $%
, where $\omega $ is the harmonic measure of $\mathcal{L}_{D}$ restricted to
the compact set $K$ given by (\ref{K}). This, in turn, implies that $\omega
|_{K}\in A_{\infty }\left( d\sigma |_{K}\right) $ and proves (\ref{enough}),
hence Theorem \ref{CRth1}, for the divergence case.

We now consider the nondivergence case. Let $\mathcal{L}_{N}=\mathbf{A}\cdot
\nabla ^{2}+\mathbf{b}\cdot \nabla $ where $\mathbf{A}$ and $\mathbf{b}$
satisfy (\ref{CarlCRa1}). Let $\Omega =\Omega _{\alpha _{0}}$, $\Phi =\Phi
_{\alpha _{0}}$, and $\rho $ be as before. Also, for $Y\in \Phi $, let $%
\mathbf{\tilde{A}}^{**}\left( Y\right) $ be as in (\ref{Ass}), and define $%
\mathcal{L}^{**}=\mathbf{A}^{**}\cdot \nabla ^{2}+\mathbf{b}\cdot \nabla $,
where 
\begin{equation*}
\mathbf{A}^{**}\left( X\right) =\left( \left( \nabla _{Y}\rho \right) \left(
\rho ^{-1}\left( X\right) \right) \right) ^{t}\mathbf{\tilde{A}}^{**}\left(
\rho ^{-1}\left( X\right) \right) \left( \nabla _{Y}\rho \right) \left( \rho
^{-1}\left( X\right) \right) \left( \det \left( \nabla _{Y}\rho \right)
\left( \rho ^{-1}\left( X\right) \right) \right) ^{-1}.
\end{equation*}
Let $\Omega _{\frac{1}{3}}=\rho \left( \Phi _{\frac{\alpha _{0}}{3}}\right) $%
, we claim that $\mathbf{A}^{**}\left( X\right) $ satisfies the following 
\begin{eqnarray}
\sup_{Z\in Q_{\frac{\delta \left( X\right) }{2\sqrt{n}}}\left( X\right) }%
\frac{\left| \mathbf{A}\left( Z\right) -\mathbf{A}^{**}\left( Z\right)
\right| ^{2}}{\delta \left( X\right) }dX &\in &\frak{C}\left( d\hat{\sigma}%
,\Omega _{\frac{1}{3}}\right) ,\qquad \text{and}  \label{AA1} \\
\sup_{Z\in Q_{\frac{\delta \left( X\right) }{2\sqrt{n}}}\left( X\right)
}\delta \left( X\right) \left| \nabla \mathbf{A}^{**}\left( Z\right) \right|
^{2}dX &\in &\frak{C}\left( d\hat{\sigma},\Omega _{\frac{1}{3}}\right)
\label{AA2}
\end{eqnarray}
where $\hat{\sigma}$ is the Lebesgue measure on $\partial \Omega _{\frac{1}{3%
}}$. Taking these properties for granted, by (\ref{AA2}), (\ref{CarlCRa1})
and Theorem \ref{KP218} applied to the operator $\mathcal{L}^{**}$, we have
that if $\omega ^{**}$ is the harmonic measure of $\mathcal{L}^{**}$ on $%
\partial \Omega _{\frac{1}{3}}$, then $\omega ^{**}\in A_{\infty }$. On the
other hand, by (\ref{AA1}) and Theorem \ref{CRth2} applied to the operators $%
\mathcal{L}$ and $\mathcal{L}^{**}$, from $\omega ^{**}\in A_{\infty }$ we
conclude that $\omega \in A_{\infty }$, where $\omega $ is the harmonic
measure of $\mathcal{L}$ on $\partial \Omega _{\frac{1}{3}}$. This finishes
the proof of Theorem \ref{CRth1} in the nondivergence case.

It only rests to establish properties (\ref{AA1}) and (\ref{AA2}).Let $Z\in
Q_{\frac{\delta \left( X\right) }{2\sqrt{n}}}\left( X\right) $ and let $%
W=\rho ^{-1}\left( Z\right) $, then from the definitions of $\mathbf{\tilde{A%
}}$ and $\mathbf{A}^{**}$ we have 
\begin{equation*}
\left| \mathbf{A}\left( Z\right) -\mathbf{A}^{**}\left( Z\right) \right|
=\left| \left( \left( \nabla _{Y}\rho \right) \left( W\right) \right)
^{t}\left[ \mathbf{\tilde{A}}\left( W\right) -\mathbf{\tilde{A}}^{**}\left(
W\right) \right] \left( \nabla _{Y}\rho \right) \left( W\right) \left( \det
\left( \nabla _{Y}\rho \right) \left( W\right) \right) ^{-1}\right| .
\end{equation*}
From (\ref{deltas}) and the fact that $\rho $ is bi-Lipschitz, and since $%
\rho ^{-1}\left( \Omega _{\frac{1}{3}}\right) =\Phi _{\frac{\alpha _{0}}{3}}$%
, 
\begin{equation*}
\frac{\left| \mathbf{A}\left( Z\right) -\mathbf{A}^{**}\left( Z\right)
\right| ^{2}}{\delta \left( Z\right) }\approx \frac{\left| \mathbf{\tilde{A}}%
\left( W\right) -\mathbf{\tilde{A}}^{**}\left( W\right) \right| ^{2}}{\tilde{%
\delta}\left( W\right) }=\frac{\left| \mathbf{\tilde{A}}^{*}\left( W\right) -%
\mathbf{\tilde{A}}^{**}\left( W\right) \right| ^{2}}{\tilde{\delta}\left(
W\right) }
\end{equation*}
Applying the proof of Lemma \ref{carlesin} to $\tilde{\delta}\left( W\right)
^{-1}\limfunc{osc}\nolimits_{\alpha \tilde{\delta}\left( W\right) }\left( 
\mathbf{\tilde{A}}^{*}-\mathbf{\tilde{A}}^{**}\right) \left( W\right) $,
from the second property in (\ref{Ezz}) it follows that for some $0<c<1$%
\begin{equation*}
\frac{\sup_{W\in Q_{c\tilde{\delta}\left( Y\right) }\left( Y\right) }\left| 
\mathbf{\tilde{A}}^{*}\left( W\right) -\mathbf{\tilde{A}}^{**}\left(
W\right) \right| ^{2}}{\tilde{\delta}\left( Y\right) }\in \frak{C}\left(
d\mu ,\Phi _{\frac{\alpha _{0}}{3}}\right) ,
\end{equation*}
where $\mu $ is the Lebesgue measure on $\partial \Phi _{\frac{\alpha _{0}}{3%
}}$; (\ref{AA1}) then follows from the fact that $\rho $ is bi-Lipschitz.
Now, by the product rule of differentiation 
\begin{eqnarray*}
\nabla \mathbf{A}^{**} &=&\nabla _{X}\left\{ \left( \nabla _{Y}\rho \right)
^{t}\mathbf{\tilde{A}}^{**}\left( \nabla _{Y}\rho \right) \det \left( \nabla
_{Y}\rho \right) ^{-1}\right\} \\
&=&\left\{ \nabla _{X}\left( \nabla _{Y}\rho \right) ^{t}\right\} \mathbf{%
\tilde{A}}^{**}\left( \nabla _{Y}\rho \right) \det \left( \nabla _{Y}\rho
\right) ^{-1} \\
&&+\left( \nabla _{Y}\rho \right) ^{t}\left\{ \nabla _{X}\mathbf{\tilde{A}}%
^{**}\left( \nabla _{Y}\rho \right) \right\} \det \left( \nabla _{Y}\rho
\right) ^{-1} \\
&&+\left( \nabla _{Y}\rho \right) ^{t}\mathbf{\tilde{A}}^{**}\left\{ \nabla
_{X}\left( \nabla _{Y}\rho \right) \right\} \det \left( \nabla _{Y}\rho
\right) ^{-1} \\
&&+\left( \nabla _{Y}\rho \right) ^{t}\mathbf{\tilde{A}}^{**}\left( \nabla
_{Y}\rho \right) \left\{ \nabla _{X}\det \left( \nabla _{Y}\rho \right)
^{-1}\right\} .
\end{eqnarray*}
Applying the chain rule in each term, we see that $\nabla \mathbf{A}^{**}$
satisfies (\ref{AA2}) because of the first property in (\ref{Ezz}), the fact
that $\tilde{\delta}\left| \nabla _{Y}^{2}\rho \right| ^{2}\in \frak{C}%
\left( d\tilde{\sigma}^{*},\Bbb{R}_{+}^{n+1}\right) $ and the boundedness of 
$\left| \nabla _{Y}\rho \right| $. 
\endproof%

\section{Proof of Theorem \ref{CRth2b}}

In the spirit of \cite{FKP} (see also \cite{CR}) we will obtain Theorem \ref
{CRth2b} as a consequence of the following perturbation result.

\begin{theorem}
\label{epsilon}Let $\mathcal{\tilde{L}}_{N,0}=\mathbf{A}\cdot \nabla ^{2}$
and $\mathcal{\tilde{L}}_{N,1}=\mathbf{A}\nabla ^{2}+\mathbf{b}\cdot \nabla $
where $\mathbf{A}=\left( A_{ij}\right) _{i,j=1}^{n+1}$ and $\mathbf{b}%
=\left( b_{j}\right) _{j=1}^{n+1}$ are bounded, measurable coefficients and $%
\mathbf{A}$ satisfy the ellipticity condition (\ref{ellipticity}). Suppose
that $\mathcal{CD}$ holds for $\mathcal{\tilde{L}}_{N,\ell }$, $\ell =0,1$.
Let $\mathcal{G}_{0}\left( X,Y\right) $ denote the Green's function for $%
\mathcal{\tilde{L}}_{N,0}$ in $\Omega $ and set $\mathcal{G}_{0}\left(
Y\right) =\mathcal{G}_{0}\left( 0,Y\right) $. There exists $\varepsilon
_{0}>0$ which depends only on $n$, $\lambda $, $\Lambda $ and $\Omega $ such
that if 
\begin{equation}
\mathcal{G}_{0}\left( X\right) \sup_{Y\in Q_{\frac{\delta \left( X\right) }{2%
\sqrt{n}}}\left( X\right) }\left| \mathbf{b}\right| ^{2}dX,\qquad X\in
\Omega ,  \label{3.1C}
\end{equation}
is a Carleson measure in $\Omega $ with respect to $d\omega _{\mathcal{%
\tilde{L}}_{N,0}}$ on $\partial \Omega $ with Carleson norm bounded by $%
\varepsilon _{0}$, i.e., 
\begin{eqnarray*}
\mathcal{G}_{0}\left( X\right) \sup_{Y\in Q_{\frac{\delta \left( X\right) }{2%
\sqrt{n}}}\left( X\right) }\left| \mathbf{b}\right| ^{2}dX &\in &\frak{C}%
\left( d\omega _{\mathcal{\tilde{L}}_{N,0}},\Omega \right) , \\
\left\| \mathcal{G}_{0}\left( X\right) \sup_{Y\in Q_{\frac{\delta \left(
X\right) }{2\sqrt{n}}}\left( X\right) }\left| \mathbf{b}\right|
^{2}dX\right\| _{\frak{C}\left( d\omega _{\mathcal{\tilde{L}}_{N,0}},\Omega
\right) } &\le &\varepsilon _{0},
\end{eqnarray*}
then $d\omega _{\mathcal{\tilde{L}}_{N,1}}\in B_{2}\left( d\omega _{\mathcal{%
\tilde{L}}_{N,0}}\right) $. Where $B_{2}\left( d\omega _{\mathcal{\tilde{L}}%
_{N,0}}\right) $ denotes the reverse H\"{o}lder class of $d\omega _{\mathcal{%
\tilde{L}}_{N,0}}$ with exponent $2$.
\end{theorem}

We defer the proof of Theorem \ref{epsilon} to the next subsection, and
prove now Theorem \ref{CRth2b}, we follow the argument in \cite{FKP}. Let $%
\triangle _{r}\left( Q\right) $ be the boundary ball $\triangle _{r}\left(
Q\right) =\left\{ P\in \partial \Omega :\left| Q-P\right| <r\right\} $, and
denote by $T_{r}\left( Q\right) $ the Carleson region in $\Omega $
associated to $\triangle _{r}\left( Q\right) $, $T_{r}\left( Q\right)
=\left\{ X\in \Omega :\left| X-Q\right| <r\right\} $. By Lemma \ref{main} we
may assume that $\mathbf{b}\left( X\right) \equiv 0$ if $\delta \left(
X\right) >r_{0}$ for some fixed (small) $r_{0}>0$. To prove Theorem \ref
{CRth2b} it is enough to show that if $\tilde{\omega}_{1}=\omega _{\mathcal{%
\tilde{L}}_{N,1}}$ with $\mathcal{\tilde{L}}_{N,1}$ as in the statement of
the theorem, then for all $Q\in \partial \Omega $, 
\begin{equation}
\tilde{\omega}_{1}|_{\triangle _{r_{0}}\left( Q\right) }\in A_{\infty }.
\label{thatsit}
\end{equation}
For $Q\in \partial \Omega ,$ $r>0$, $\alpha >0$, let $\Gamma _{\alpha
,r}\left( Q\right) $ be a nontangential cone of fixed aperture $\alpha $ and
height $r$, i.e. 
\begin{equation*}
\Gamma _{\alpha ,r}\left( Q\right) =\left\{ X\in \Omega :\left| X-Q\right|
<\left( 1+\alpha \right) \delta \left( X\right) <\left( 1+\alpha \right)
r\right\} .
\end{equation*}
For a fixed $\alpha _{0}>0$ to be determined later, let $\mathcal{E}%
_{r}\left( Q\right) $ be given by 
\begin{equation*}
\mathcal{E}_{\mathbf{b},r}\left( Q\right) =\left\{ \int_{\Gamma _{\alpha
_{0},r}\left( Q\right) }\delta \left( X\right) ^{1-n}\sup_{Y\in Q_{\frac{%
\delta \left( X\right) }{2\sqrt{n}}}\left( X\right) }\left| \mathbf{b}%
\right| ^{2}dX\right\} ^{\frac{1}{2}},\qquad Q\in \partial \Omega .
\end{equation*}
Fix $\alpha _{0}$ , $r_{0}$, such that $\Gamma _{\alpha _{0},r_{0}}\left(
Q\right) \subset T_{2r_{0}}\left( Q\right) $ for all $Q\in \partial \Omega $%
. Then, letting $\sigma $ be the Lebesgue measure on $\partial \Omega $, by
Fubini's theorem, the hypothesis $\delta \left( X\right) \sup_{Y\in Q_{\frac{%
\delta \left( X\right) }{2\sqrt{n}}}\left( X\right) }\left| \mathbf{b}\left(
Y\right) \right| ^{2}dX\in \frak{C}\left( d\sigma ,\Omega \right) $ and the
doubling property of $d\sigma $, we have 
\begin{equation*}
\frac{1}{\sigma \left( \triangle _{r_{0}}\left( Q\right) \right) }%
\int_{\triangle _{r_{0}}\left( Q\right) }\mathcal{E}_{\mathbf{b}%
,r_{0}}\left( P\right) ^{2}d\sigma \left( P\right) \le C\frac{1}{\sigma
\left( \triangle _{r_{0}}\left( Q\right) \right) }\int_{T_{2r_{0}}\left(
Q\right) }\delta \left( X\right) \sup_{Y\in Q_{\frac{\delta \left( X\right) 
}{2\sqrt{n}}}\left( X\right) }\left| \mathbf{b}\left( Y\right) \right|
^{2}dX\le C.
\end{equation*}
That is, the average in $\triangle _{r_{0}}\left( Q\right) $ of $\mathcal{E}%
_{\mathbf{b},r_{0}}^{2}$ is bounded. Hence, there exists a closed set $%
F\subset \triangle _{r_{0}}\left( Q\right) $\thinspace such that $\sigma
\left( F\right) >\frac{1}{2}\sigma \left( \triangle _{r0}\left( Q\right)
\right) $ and $\mathcal{E}_{\mathbf{b},r_{0}}\left( P\right) \le C$ for all $%
P\in F$. Let $\Omega _{F}$ be a ``saw-tooth'' region above $F$ in $\Omega $
as given in definition \ref{serrucho}, and let $\mathbf{b}^{*}\left(
X\right) $ be given by 
\begin{equation*}
\mathbf{b}^{*}\left( X\right) =\left\{ 
\begin{array}{lll}
\mathbf{b}\left( X\right) &  & X\in \Omega _{F} \\ 
\mathbf{0} &  & X\in \Omega \backslash \Omega _{F}
\end{array}
\right. .
\end{equation*}
The drift $\mathbf{b}^{*}$ so defined satisfies $\mathcal{E}_{\mathbf{b}%
^{*},r_{0}}\left( P\right) \le C_{0}$ for \emph{all} $P\in \triangle
_{r_{0}}\left( Q\right) $. We claim that if $C_{0}\,$is small enough then
the operators $\mathcal{\mathcal{L}}_{0}=\mathcal{\tilde{L}}_{N}=\mathbf{A}%
\cdot \nabla ^{2}$ and $\mathcal{L}_{1}=\mathbf{A}\cdot \nabla ^{2}+\mathbf{b%
}^{*}\cdot \nabla $ satisfy the hypotheses of Theorem \ref{epsilon} in $\Phi
=T_{3r_{0}}\left( Q\right) $. Indeed, since $\mathbf{b}^{*}\equiv 0$ in $%
T_{3r_{0}\left( Q\right) }\backslash T_{2r_{0}}\left( Q\right) $, we only
need to check the Carleson measure condition (\ref{3.1C}) near $\triangle
_{2r_{0}}\left( Q\right) $. More precisely, we will show that for all $%
s<r_{0}/2$ and $P\in \triangle _{2r_{0}}\left( Q\right) $, 
\begin{equation}
\int_{T_{s}\left( P\right) }\mathcal{G}_{0}\left( X\right) \sup_{Y\in Q_{%
\frac{\delta \left( X\right) }{2\sqrt{n}}}\left( X\right) }\left| \mathbf{b}%
\right| ^{2}dX\le \varepsilon _{0}\omega _{0}\left( \triangle _{s}\left(
P\right) \right) .  \label{3.1CC}
\end{equation}
where $\mathcal{G}_{0}\left( X,Y\right) $ is the Green's function for $%
\mathcal{L}_{0}$ in $\Phi $, $\mathcal{G}_{0}\left( Y\right) =\mathcal{G}%
_{0}\left( X_{0},Y\right) $ (where $X_{0}\,$is the \emph{center} of $\Phi $)
and $\omega _{0}=\omega _{\mathcal{\mathcal{L}}_{0},\Phi }^{X_{0}}$ in $%
\partial \Phi $. As in \cite{CR} (see Lemmas 2.8 and 2.14 there), we have 
\begin{equation}
\frac{1}{\omega _{0}\left( \triangle _{X}\right) }\frac{\mathcal{G}%
_{0}\left( X\right) }{\delta \left( X\right) ^{2}}\approx \frac{\mathcal{G}%
\left( X\right) }{\int_{Q_{\frac{\delta \left( X\right) }{2\sqrt{n}}}}%
\mathcal{G}\left( Y\right) dY},  \label{2.8}
\end{equation}
where $\triangle _{X}=\triangle _{\delta \left( X\right) }\left( Q\right) $
for some $Q\in \partial \Phi $\thinspace such that $\left| X-Q\right|
=\delta \left( X\right) ;$ and $\mathcal{G}\left( X\right) =\mathcal{G}%
\left( \overline{X},X\right) $, with $\mathcal{G}\left( Z,X\right) $ the
Green's function for $\mathcal{L}_{0}$ in a fixed domain $\mathcal{T}\supset
\supset \Omega $ and $\overline{X}\in \mathcal{T}\backslash \overline{\Omega 
}$ is a fixed point away from $\Omega $. From (\ref{2.8}), writing $\mathcal{%
G}\left( Q_{\frac{\delta \left( X\right) }{2\sqrt{n}}}\right) =\int_{Q_{%
\frac{\delta \left( X\right) }{2\sqrt{n}}}}\mathcal{G}\left( Y\right) dY$
and proceeding as in the proof of (5.1) in \cite{CR} we have 
\begin{eqnarray*}
\int_{T_{s}\left( P\right) }\mathcal{G}_{0}\left( X\right) \sup_{Y\in Q_{%
\frac{\delta \left( X\right) }{2\sqrt{n}}}\left( X\right) }\left| \mathbf{b}%
^{*}\right| ^{2}dX &\le &\int_{T_{s}\left( P\right) }\omega _{0}\left(
\triangle _{X}\right) \frac{\delta \left( X\right) ^{2}\mathcal{G}\left(
X\right) }{\mathcal{G}\left( Q_{\frac{\delta \left( X\right) }{2\sqrt{n}}%
}\right) }\sup_{Y\in Q_{\frac{\delta \left( X\right) }{2\sqrt{n}}}\left(
X\right) }\left| \mathbf{b}\right| ^{2}dX \\
&\le &C\int_{\triangle _{s}\left( P\right) }\mathcal{E}_{\mathbf{b}%
^{*},r}^{2}\left( P\right) d\omega _{0}\left( P\right) \\
&\le &CC_{0}^{2}\omega _{0}\left( \triangle _{s}\left( P\right) \right) \le
\varepsilon _{0}\omega _{0}\left( \triangle _{s}\left( P\right) \right)
\end{eqnarray*}
if $C_{0}$ is small enough.Thus, $\mathcal{\mathcal{L}}_{0}$ and $\mathcal{L}%
_{1}$ satisfy the hypotheses of Theorem \ref{epsilon} in $\Phi $ and hence $%
\omega _{1}=\omega _{\mathcal{L}_{1},\Phi }^{X_{0}}\in B_{2}\left( d\omega
_{0},\partial \Phi \right) $. By property (2) in Lemma \ref{pro} we have
that for $E\subset \triangle _{2r_{0}}\left( Q\right) $, 
\begin{equation*}
\omega _{0}\left( E\right) \approx \frac{\omega _{\mathcal{L}_{0},\Omega
}\left( E\right) }{\omega _{\mathcal{L}_{0},\Omega }\left( \triangle \right) 
}.
\end{equation*}
This together with the hypothesis from Theorem \ref{CRth2b} that $\omega _{%
\mathcal{L}_{0},\Omega }\in A_{\infty }$, implies $\omega _{0}|_{\triangle
_{3r_{0}}\left( Q\right) }\in A_{\infty }$. From $\omega _{1}\in B_{2}\left(
d\omega _{0},\partial \Phi \right) $ we conclude that $\omega
_{1}|_{\triangle _{2r_{0}}\left( Q\right) }\in A_{\infty }$. Hence, for some
constants $0<\alpha _{0},c_{0}$, we have 
\begin{equation}
\left( \omega _{1}\left( F\right) \right) ^{\alpha _{0}}\approx \left( \frac{%
\omega _{1}\left( F\right) }{\omega _{1}\left( \triangle _{r_{0}}\left(
Q\right) \right) }\right) ^{\alpha _{0}}\ge c_{0}\frac{\sigma \left(
F\right) }{\sigma \left( \triangle _{r_{0}}\left( Q\right) \right) }\ge 
\frac{c_{0}}{2},  \label{that}
\end{equation}
where we applied property (1) of Lemma \ref{pro} and we used that $\sigma
\left( F\right) >\frac{1}{2}\sigma \left( \triangle _{r0}\left( Q\right)
\right) $.

Now, let $\nu $ be the harmonic measure of $\mathcal{L}_{1}$ in $\Omega _{F}$%
. By Lemma \ref{mmain} we have that for some $0<\theta <1$, 
\begin{equation}
\frac{\omega _{1}\left( F\right) }{\omega _{1}\left( \triangle
_{r_{0}}\left( Q\right) \right) }<\nu \left( F\right) ^{\theta }.
\label{teta}
\end{equation}
On the other hand, since $\mathbf{b}^{*}$ coincides with $\mathbf{b}$ in $%
\Omega _{F}$, any solution $u$ of $\mathcal{\tilde{L}}_{1}u=0$ in $\Phi \,$%
is a solution of $\mathcal{\tilde{L}}_{1}u=0$ in $\Omega _{F}\subset \Phi $.
The boundary maximum principle implies that for all $F\subset \triangle
_{r_{0}}\left( Q\right) $, $\nu \left( F\right) \le \omega _{1}^{*}\left(
F\right) $, where $\omega _{1}^{*}=\omega _{\mathcal{\tilde{L}}_{1},\Phi }$.
From (\ref{teta}) and (\ref{that}) then we obtain 
\begin{equation*}
\omega _{1}^{*}\left( F\right) >c_{1}>0.
\end{equation*}
By Property (2) in Lemma \ref{pro} and the maximum principle, we have 
\begin{equation*}
\frac{\tilde{\omega}_{1}\left( F\right) }{\tilde{\omega}_{1}\left( \triangle
_{r_{0}}\left( Q\right) \right) }\ge \omega _{1}^{*}\left( F\right) >c_{1}>0.
\end{equation*}
Therefore, whenever $\frac{\sigma \left( F\right) }{\sigma \left( \triangle
_{r_{0}}\left( Q\right) \right) }>\frac{1}{2}$ then$\frac{\tilde{\omega}%
_{1}\left( F\right) }{\tilde{\omega}_{1}\left( \triangle _{r_{0}}\left(
Q\right) \right) }>c_{1}$ . This shows that (\ref{thatsit}) holds. 
\endproof%

\subsection{Proof of Theorem \ref{epsilon}}

The proof of this results closely follows the steps in \cite{CR}. We will
sketch the main steps and refer the reader to \cite{CR} for the technical
details omitted here. First, by standard arguments the problem is reduced to
treating the case in which $\Omega $\thinspace is the unit ball $%
B=B_{1}\left( 0\right) $ (this is justified as far as the methods are \emph{%
preserved} under bi-Lipschitz transformation). For simplicity, we will write 
$\omega _{0}=\omega _{\mathcal{\tilde{L}}_{N,0}}$ and $\omega _{1}=$ $\omega
_{\mathcal{\tilde{L}}_{N,1}}$. To see that $d\omega _{1}\in B_{2}\left(
d\omega _{0}\right) $ it is equivalent to prove that if $u_{1}$ is a
solution of the Dirichlet problem 
\begin{equation*}
\left\{ 
\begin{array}{lll}
\mathcal{\tilde{L}}_{N,1}u_{1} & = & 0\qquad \text{in }B \\ 
u_{1} & = & g\qquad \text{on }\partial B,
\end{array}
\right.
\end{equation*}
where $g$ is continuous in $\partial B$, then 
\begin{equation}
\left\| Nu_{1}\right\| _{L^{2}\left( \partial B,d\omega _{0}\right) }\le
C\left\| g\right\| _{L^{2}\left( \partial B,d\omega _{0}\right) },
\label{Neps}
\end{equation}
where $Nu$ is the nontangential maximal function (with some fixed aperture $%
\alpha >0$) of $u$. We let $u_{0}$ be the solution to 
\begin{equation*}
\left\{ 
\begin{array}{lll}
\mathcal{\tilde{L}}_{N,0}u_{0} & = & 0\qquad \text{in }B \\ 
u_{0} & = & g\qquad \text{on }\partial B.
\end{array}
\right.
\end{equation*}
Then $u_{1}-u_{0}=0$ on $\partial B$ and we have the representation 
\begin{equation}
u_{1}\left( X\right) =u_{0}\left( X\right) -\int_{B}\mathcal{G}_{0}\left(
X,Y\right) \mathcal{\tilde{L}}_{N,0}u_{1}dY=u_{0}\left( X\right) -F\left(
X\right) .  \label{F}
\end{equation}
Then, (\ref{Neps}) follows as in \cite{CR} from the following two lemmas.

\begin{lemma}
\label{3.2}Let $\mathcal{G}\left( X,Y\right) $ denote the Green's function
for $\mathcal{\tilde{L}}_{N,0}$ in $B_{10}\left( 0\right) $ and let $%
\mathcal{G}\left( Y\right) =\mathcal{G}\left( \bar{X},Y\right) $ where $\bar{%
X}$ is some fixed point in $\Bbb{R}^{n+1}$ such that $\left| \bar{X}\right|
=5$. For $Y\in \mathcal{B}$ let $B_{0}\left( Y\right) $ and $B\left(
Y\right) $ denote the Euclidean balls centered at $Y\,$of radii $\frac{%
\delta \left( Y\right) }{6}$ and $\frac{\delta \left( Y\right) }{2}$
respectively (in this case $\delta \left( Y\right) =\limfunc{dist}\left(
Y,\partial B\right) =1-\left| Y\right| $). Under the hypotheses of Theorem 
\ref{epsilon}, we have that $F$ as in (\ref{F}) satisfies 
\begin{equation*}
N^{0}F\left( Q\right) =\sup_{X\in \Gamma \left( Q\right) }\left\{
\int_{B_{0}\left( X\right) }F^{2}\left( Y\right) \frac{\mathcal{G}\left(
Y\right) }{\mathcal{G}\left( B\left( Y\right) \right) }dY\right\} ^{\frac{1}{%
2}}\le C\varepsilon _{0}M_{\omega _{0}}\left( Su_{1}\right) \left( Q\right) ,
\end{equation*}
where 
\begin{equation*}
M_{\omega _{0}}f\left( Q\right) =\sup_{r>0}\frac{1}{\omega _{0}\left(
\triangle _{r}\left( Q\right) \right) }\int_{\triangle _{r}\left( Q\right)
}\left| f\left( P\right) \right| d\omega _{0}
\end{equation*}
is the Hardy-Littlewood maximal function of $f$ with respect to the measure $%
\omega _{0\text{,}}$, and $Su_{1}$ is the area function of $u_{1}$.
\end{lemma}

\begin{lemma}
\label{3.3}Under the hypotheses of Lemma \ref{3.2}, 
\begin{equation*}
\int_{\partial B}S^{2}u_{1}d\omega _{0}\le C\int_{\partial B}\left(
Nu_{1}\right) ^{2}d\omega _{0}.
\end{equation*}
\end{lemma}

Indeed, by Lemma 2.21 in \cite{CR} it follows that 
\begin{eqnarray*}
\int_{\partial B}\left( Nu_{1}\right) ^{2}d\omega _{0} &\le &C\int_{\partial
B}\left( N^{0}u_{1}\right) ^{2}d\omega _{0} \\
&\le &C\int_{\partial B}\left( N^{0}u_{0}\right) ^{2}d\omega
_{0}+C\int_{\partial B}\left( N^{0}F\right) ^{2}d\omega _{0}
\end{eqnarray*}
where $N^{0}u_{1}\,$is as in Lemma \ref{3.2}. Given that Lemmas \ref{3.2}
and \ref{3.3} hold, we have 
\begin{eqnarray*}
\int_{\partial B}\left( Nu_{1}\right) ^{2}d\omega _{0} &\le &C\int_{\partial
B}\left( Nu_{0}\right) ^{2}d\omega _{0}+C\varepsilon _{0}\int_{\partial
B}M_{\omega _{0}}\left( Su_{1}\right) ^{2}\left( Q\right) d\omega _{0} \\
&\le &C\int_{\partial B}g^{2}d\omega _{0}+C\varepsilon _{0}\int_{\partial
B}\left( Su_{1}\right) ^{2}\left( Q\right) d\omega _{0} \\
&\le &C\int_{\partial B}g^{2}d\omega _{0}+C\varepsilon _{0}\int_{\partial
B}\left( Nu_{1}\right) ^{2}\left( Q\right) d\omega _{0}
\end{eqnarray*}
and the last term on the right can be absorbed into the left if $\varepsilon
_{0}$ is small enough. This proves (\ref{Neps}) and hence Theorem \ref
{epsilon}.%
\endproof%

We will only write in some detail the proof of Lemma \ref{3.2}. Given the
big overlap with the methods in \cite{CR} this will suffice to indicate the
proof of Lemma \ref{3.3}, which is almost identical to the proof of Lemma
3.3 in \cite{CR}.

\proof of Lemma~\ref{3.2}%
Fix $Q_{0}\in \partial B$ and $X_{0}\in \Gamma \left( Q_{0}\right) $. Let $%
B_{0}=B_{\frac{\delta _{0}}{6}}\left( X_{0}\right) $ and $KB_{0}=B_{\frac{%
K\delta _{0}}{6}}\left( X_{0}\right) $ where $\delta _{0}=\delta \left(
X_{0}\right) $ and $K>0$. Let $\mathcal{\tilde{G}}\left( X,Y\right) $ be the
Green's function for $\mathcal{\tilde{L}}_{N,0}$ on $3B_{0}$, set 
\begin{equation*}
\begin{array}{lll}
F_{1}\left( X\right) & = & \int_{2B_{0}}\mathcal{\tilde{G}}\left( X,Y\right) 
\mathcal{\tilde{L}}_{N,0}u_{1}\left( Y\right) dY, \\ 
&  &  \\ 
F_{2}\left( X\right) & = & \int_{2B_{0}}\left[ \mathcal{G}_{0}\left(
X,Y\right) -\mathcal{\tilde{G}}\left( X,Y\right) \right] \mathcal{\tilde{L}}%
_{N,0}u_{1}\left( Y\right) dY, \\ 
&  &  \\ 
F_{3}\left( X\right) & = & \int_{B\backslash 2B_{0}}\mathcal{G}_{0}\left(
X,Y\right) \mathcal{\tilde{L}}_{N,0}u_{1}\left( Y\right) dY.
\end{array}
\end{equation*}

So that $F$ in (\ref{F}) is given by $F\left( X\right) =F_{1}\left( X\right)
+F_{2}\left( X\right) +F_{3}\left( X\right) $, and proving Lemma \ref{3.2}
is reduced to proving that 
\begin{equation*}
\int_{_{B_{0}}}F_{i}^{2}\left( Y\right) \frac{\mathcal{G}\left( Y\right) }{%
\mathcal{G}\left( B\left( Y\right) \right) }dY\le C\varepsilon
_{0}^{2}M_{\omega _{0}}^{2}\left( Su_{1}\right) \left( Q_{0}\right) ,\qquad 
\text{for }i=1,2,3.
\end{equation*}

We will only prove this in some detail for $i=1$. Even though $i=1$ is
allegedly the simplest case of the three, its proof captures the significant
differences with the proof of Lemma 3.2 in \cite{CR} (the analog to Lemma 
\ref{3.2} here), so the other two cases follow in a similar manner as in 
\cite{CR}.

Let $\beta \left( X\right) =\sup_{\left| Z-X\right| \le \frac{\delta \left(
X\right) }{2}}\left| \mathbf{b}\left( Z\right) \right| $, then, given $Y\in
B_{0}$, we have 
\begin{eqnarray}
&&\left| \mathbf{b}\left( Y\right) \right|  \notag \\
&\le &\frac{C}{\left| 3B_{0}\left( Y\right) \right| }\int_{3B_{0}\left(
Y\right) }\beta \left( X\right) dX  \notag \\
&\le &\frac{C}{\left| 3B_{0}\left( Y\right) \right| }\left\{
\int_{3B_{0}\left( Y\right) }\beta \left( X\right) ^{2}\frac{\mathcal{G}%
\left( X\right) }{\mathcal{G}\left( B\left( X\right) \right) }dX\right\} ^{%
\frac{1}{2}}\left\{ \int_{3B_{0}\left( Y\right) }\frac{\mathcal{G}\left(
B\left( X\right) \right) }{\mathcal{G}\left( X\right) }dX\right\} ^{\frac{1}{%
2}}.  \label{pupo}
\end{eqnarray}
Using (\ref{2.8}) on the right side of (\ref{pupo}), applying the Carleson
measure property of $\mathcal{G}_{0}\left( X\right) \beta \left( X\right)
^{2}$ and the doubling property of $\omega _{0}$, we obtain 
\begin{eqnarray}
&&\left| \mathbf{b}\left( Y\right) \right|  \notag \\
&\le &\frac{C}{\left| 3B_{0}\left( Y\right) \right| }\left\{ \frac{1}{\omega
_{0}\left( \triangle _{Y}\right) }\int_{3B_{0}\left( X\right) }\beta \left(
X\right) ^{2}\frac{\mathcal{G}_{0}\left( X\right) }{\delta \left( X\right)
^{2}}dX\right\} ^{\frac{1}{2}}\left\{ \int_{3B_{0}\left( Y\right) }\frac{%
\mathcal{G}\left( B\left( X\right) \right) }{\mathcal{G}\left( X\right) }%
dX\right\} ^{\frac{1}{2}}  \notag \\
&\le &\frac{C}{\left| 3B_{0}\left( Y\right) \right| }\left\{ \frac{1}{\omega
_{0}\left( \triangle _{Y}\right) }\int_{3B_{0}\left( Y\right) }\beta \left(
X\right) ^{2}\frac{\mathcal{G}_{0}\left( X\right) }{\delta \left( X\right)
^{2}}dX\right\} ^{\frac{1}{2}}\left\{ \int_{3B_{0}\left( Y\right) }\frac{%
\mathcal{G}\left( B\left( X\right) \right) }{\mathcal{G}\left( X\right) }%
dX\right\} ^{\frac{1}{2}}  \notag \\
&\le &\frac{C\varepsilon _{0}}{\left| 3B_{0}\left( Y\right) \right| \delta
_{0}}\left\{ \int_{3B_{0}\left( Y\right) }\frac{\mathcal{G}\left( B\left(
X\right) \right) }{\mathcal{G}\left( X\right) }dX\right\} ^{\frac{1}{2}}\le 
\frac{C\varepsilon _{0}}{\delta _{0}},  \label{bb}
\end{eqnarray}
where the last inequality follows as (3.6) in \cite{CR}. Note that $F_{1}$
satisfies $\mathcal{\tilde{L}}_{N,0}F_{1}=\chi \left( 2B_{0}\right) \mathcal{%
\tilde{L}}_{N,0}u_{1}$ in $3B_{0}$, where $\chi \left( 2B_{0}\right) $
denotes the characteristic function of $2B_{0}$. Hence, from the (weighted)
a priori estimates for solutions (Theorem 2.5 in \cite{CR}), 
\begin{equation}
\left\{ \int_{_{3B_{0}}}\left| \nabla F_{1}\left( Y\right) \right| ^{2}\frac{%
\mathcal{G}\left( Y\right) }{\mathcal{G}\left( B\left( Y\right) \right) }%
dY\right\} ^{\frac{1}{2}}\le C\delta _{0}\left\{ \int_{_{3B_{0}}}\left| 
\mathcal{\tilde{L}}_{N,0}F_{1}\left( Y\right) \right| ^{2}\frac{\mathcal{G}%
\left( Y\right) }{\mathcal{G}\left( B\left( Y\right) \right) }dY\right\} ^{%
\frac{1}{2}}.  \label{pu}
\end{equation}
On the other hand, by the weighted Poincar\'{e} inequality (Theorem 1.2 in 
\cite{FKS}), 
\begin{equation*}
\left\{ \int_{_{3B_{0}}}F_{1}^{2}\left( Y\right) \frac{\mathcal{G}\left(
Y\right) }{\mathcal{G}\left( B\left( Y\right) \right) }dY\right\} ^{\frac{1}{%
2}}\le C\delta _{0}\left\{ \int_{3B_{0}}\left| \nabla F_{1}\left( Y\right)
\right| ^{2}\frac{\mathcal{G}\left( Y\right) }{\mathcal{G}\left( B\left(
Y\right) \right) }dY\right\} ^{\frac{1}{2}}.
\end{equation*}
Combining this with (\ref{pu}), using 
\begin{equation*}
\mathcal{\tilde{L}}_{N,0}F_{1}=\chi \left( 2B_{0}\right) \left\{ \mathcal{%
\tilde{L}}_{N,0}-\mathcal{\tilde{L}}_{N,1}\right\} u_{1}=\chi \left(
2B_{0}\right) \mathbf{b}\cdot \nabla u_{1},
\end{equation*}
and (\ref{bb}), we get 
\begin{eqnarray}
\left\{ \int_{_{B_{0}}}F_{1}^{2}\left( Y\right) \frac{\mathcal{G}\left(
Y\right) }{\mathcal{G}\left( B\left( Y\right) \right) }dY\right\} ^{\frac{1}{%
2}} &\le &C\delta _{0}^{2}\left\{ \int_{_{3B_{0}}}\left| \mathcal{\tilde{L}}%
_{N,0}F_{1}\left( Y\right) \right| ^{2}\frac{\mathcal{G}\left( Y\right) }{%
\mathcal{G}\left( B\left( Y\right) \right) }dY\right\} ^{\frac{1}{2}}  \notag
\\
&\le &C\left\{ \int_{_{2B_{0}}}\delta \left( Y\right) ^{4}\left| \mathbf{b}%
\right| ^{2}\left| \nabla u_{1}\right| ^{2}\frac{\mathcal{G}\left( Y\right) 
}{\mathcal{G}\left( B\left( Y\right) \right) }dY\right\} ^{\frac{1}{2}} 
\notag \\
&\le &C\varepsilon _{0}\left\{ \int_{_{2B_{0}}}\delta \left( Y\right)
^{2}\left| \nabla u_{1}\right| ^{2}\frac{\mathcal{G}\left( Y\right) }{%
\mathcal{G}\left( B\left( Y\right) \right) }dY\right\} ^{\frac{1}{2}}  \notag
\\
&\le &C\varepsilon _{0}Su_{1}\left( Q_{0}\right) .  \label{F1}
\end{eqnarray}
The rest of the proof proceeds as in \cite{CR}, to obtain 
\begin{eqnarray}
\left\{ \int_{_{B_{0}}}F_{2}^{2}\left( Y\right) \frac{\mathcal{G}\left(
Y\right) }{\mathcal{G}\left( B\left( Y\right) \right) }dY\right\} ^{\frac{1}{%
2}} &\le &C\varepsilon _{0}Su_{1}\left( Q_{0}\right) \qquad \text{and}
\label{F2} \\
\left\{ \int_{_{B_{0}}}F_{3}^{2}\left( Y\right) \frac{\mathcal{G}\left(
Y\right) }{\mathcal{G}\left( B\left( Y\right) \right) }dY\right\} ^{\frac{1}{%
2}} &\le &C\varepsilon _{0}M_{\omega _{0}}\left( Su_{1}\right) \left(
Q_{0}\right)  \label{F3}
\end{eqnarray}
respectively. Since $F\left( Y\right) =F_{1}\left( Y\right) +F_{2}\left(
Y\right) +F_{3}\left( Y\right) $, we have 
\begin{equation*}
\int_{B_{0}\left( X_{0}\right) }F^{2}\left( Y\right) \frac{\mathcal{G}\left(
Y\right) }{\mathcal{G}\left( B\left( Y\right) \right) }dY\le
C\sum_{i=1}^{3}\int_{B_{0}\left( X_{0}\right) }F_{i}^{2}\left( Y\right) 
\frac{\mathcal{G}\left( Y\right) }{\mathcal{G}\left( B\left( Y\right)
\right) }dY.
\end{equation*}
Lemma \ref{3.2} then follows from (\ref{F1})--(\ref{F3}) by taking supremum
over all $X_{0}\in \Gamma \left( Q_{0}\right) $. 
\endproof%


\begin{thebibliography}{99}
\bibitem{BA}  A. Beurling and L. Ahlfors, \emph{The boundary correspondence
under quasi-conformal mapping}, Acta. Math. \textbf{96} (1956), 125--142.

\bibitem{BAR}  T. Barcel\'{o}, \emph{A comparison and Fatou theorem for a
class of nondivergence elliptic equations with singular lower order terms}.
Indiana Univ. Math. J. \textbf{43} (1994), 1, 1--24

\bibitem{CFL}  F. Chiarenza, M. Frasca and P. Longo, \emph{$W^{2,p}$%
-solvability for the Dirichlet problem for nondivergence elliptic equations
with VMO coefficients,} Trans. AMS \textbf{336}, $2$ (1993), 841--853.

\bibitem{DJK}  B. Dahlberg, D. Jerison and C. Kenig, \emph{Area integral
estimates for elliptic differential operators with nonsmooth coefficients,}%
Arkiv. Mat \textbf{22} (1984), 97-108

\bibitem{EK}  L. Escauriaza and C. Kenig, \emph{Area integral estimates for
for solutions and normalized adjoint solutions to nondivergence form
elliptic equations}, Arkiv. Mat. 31 (1993), 275-296.

\bibitem{FKP}  R. Fefferman, C. Kenig C. and J. Pipher, \emph{The theory of
weights and the Dirichlet problem for elliptic equations}, Annals of Math. 
\textbf{134} (1991), 65-124.

\bibitem{FKS}  Fabes, Eugene B.; Kenig, Carlos E.; Serapioni, Raul P. \emph{%
The local regularity of solutions of degenerate elliptic equations.} Comm.
Partial Differential Equations 7 (1982), no. 1, 77--116.

\bibitem{HL}  S. Hoffman and J. Lewis, \emph{The Dirichlet problem for
parabolic operators with singular drift terms}, Mem. Amer. Math. Soc. 151, 
\textbf{719 }(2001)

\bibitem{KKPT}  C. Kenig, H. Koch, J. Pipher and T. Toro, \emph{A New
Approach to Absolute Continuity of Elliptic Measure, with Applications to
Non-symmetric Equations}, Advances in Mathematics \textbf{153}, (2000)
231--298.

\bibitem{KP}  C. Kenig and J. Pipher, \emph{The Dirichlet problem for
elliptic equations with drift terms}, Publ. Mat. \textbf{45} (2001), no. 1,
199--217.

\bibitem{Ki}  T. Kilpel\"{a}inen, \emph{Weighted Sobolev spaces and capacity}%
, Ann. Acad. Sci. Fenn. Ser. A I. Math. \textbf{19} (1994), 95-113.

\bibitem{LV}  O. Lehto and K. Virtanen, \emph{Quasiconformal Mappings in the
Plane}, second edition, Springer-Verlag 1973.

\bibitem{MMo}  L. Modica and S. Mortola, \emph{Construction of a singular
elliptic-harmonic-measure,} Manuscrita Math. \textbf{33} (1980), 81--98.

\bibitem{Mo}  J. Moser, \emph{On Harnack's theorem for elliptic differential
equations}. Comm. Pure and Appl. Math. \textbf{14} (1961), 577-591.

\bibitem{Mu}  B. Muckenhoupt, \emph{The equivalence of two conditions for
weight functions}, Studia Math. \textbf{49} (1974), 101--106.

\bibitem{CRt}  C. Rios, \emph{Sufficient Conditions for the absolute
continuity of the nondivergence harmonic measure}, Ph.D. Thesis, University
of Minnesota, Minneapolis, (2001).

\bibitem{CR}  C. Rios, \emph{The }$L^{p}$\emph{\ Dirichlet Problem and
Nondivergence Harmonic Measure}, Trans. AMS \textbf{355}, 2 (2003), 665--687.

\bibitem{S}  M. V. Safonov, \emph{Nonuniqueness for second order elliptic
equations with measurable coefficients}, SIAM J. Math. Anal. \textbf{30}
(1999), pp. 879--895
\end{thebibliography}
\end{document}